\documentclass{amsart}
\usepackage{amssymb}
\usepackage{amsmath}
\usepackage{amsthm}
\usepackage[all]{xy}
\xyoption{2cell}

\allowdisplaybreaks                          
\newcommand{\F}{{\mathcal F}}
\newcommand{\Ol}{{\mathcal O}}
\newcommand{\f}{\varphi}
\newcommand{\ft}{\varphi_\theta}
\newcommand{\fot}{\varphi_{\overline{\theta}}}
\newcommand{\pst}{\psi_\vartheta}
\newcommand{\psot}{\psi_{\overline{\vartheta}}}

\newcommand{\E}{{\mathcal E}}

\newcommand{\G}{{\mathcal G}}

\newcommand{\pu}{{\mathbb P^1}}
\newcommand{\proj}{\mathbb P}
\newcommand{\quadr}{\mathbb Q}

\newcommand{\pt}{{\mathbb P^3}}

\newcommand{\pd}{{\mathbb P^2}}
\newcommand{\qt}{{\mathbb Q^3}}

\newcommand{\lra}{{\longrightarrow}}
\newcommand{\shse}[3]{0 ~\lra ~#1~ \lra ~#2~ \lra ~#3~ \lra~ 0}
\DeclareMathOperator{\loc}{\mathrm{Locus}}
\DeclareMathOperator{\cloc}{\mathrm{ChLocus}}
\newcommand{\ratcurves}{\textrm{Ratcurves}^n(X)}

\newcommand{\om}{\textrm{Hom}}
\newcommand{\Aut}{\textrm{Aut}}
\newcommand{\Univ}{\textrm{Univ}}
\newcommand{\cone}{\textrm{NE}}
\DeclareMathOperator{\pic}{Pic}

\DeclareMathOperator{\Exc}{Exc}
\DeclareMathOperator{\rk}{\mathrm rk}

\newcommand{\oth}{{\overline{\theta}}}
\newcommand{\ovth}{{\overline{\vartheta}}}
\newcommand{\ba}{{B\v anic\v a~}}
\begin{document}

\newtheorem{theorem}{Theorem}[section]
\newtheorem{lemma}[theorem]{Lemma}
\newtheorem{proposition}[theorem]{Proposition}
\newtheorem{corollary}[theorem]{Corollary}
\newtheorem{example}[theorem]{Example}
\theoremstyle{definition}
\newtheorem{definition}[theorem]{Definition}
\newtheorem{statement}[theorem]{}
\newtheorem{xca}[theorem]{Exercise}
\theoremstyle{remark}
\newtheorem{remark}[theorem]{Remark}

\renewcommand{\theequation}{{\arabic{section}.\arabic{theorem}.\arabic{equation}}
}
\author{Carla Novelli}

\author{Gianluca Occhetta}

\address{Dipartimento di Matematica, via Sommarive 14, I-38050 Povo (TN)}
\email{novelli@science.unitn.it}
\email{occhetta@science.unitn.it}

\title{Ruled Fano fivefolds of index two}

\begin{abstract}
We classify Fano fivefolds of index two which are projectivization
of rank two vector bundles over four dimensional manifolds.
\end{abstract}

\maketitle


\section{Introduction}

A smooth complex projective variety $X$ is called {\sf Fano} if its
anticanonical bundle $-K_X$ is ample; the {\sf index} of $X$, $r_X$, 
is the largest natural number $m$ such that $-K_X=mH$ for some (ample) 
divisor $H$ on $X$, while the {\sf pseudoindex}, $i_X$, is the minimum 
anticanonical degree of rational curves on $X$.\\
Since $X$ is smooth, $\pic(X)$ is torsion free, and therefore the divisor
$L$ satisfying $-K_X = r_XL$ is uniquely determined and called the 
{\sf fundamental divisor} of $X$.\\
By a theorem of Kobayashi and Ochiai \cite{Koba}, $r_X \ge \dim X +1$ if and only if 
$(X,L) \simeq (\proj^{\dim X}, \Ol_{\proj}(1))$, and $r_X= \dim X$ if and only if 
$(X,L) \simeq (\quadr^{\dim X}, \Ol_{\quadr}(1))$.\\
Fano manifolds of index $\dim X-1$ and $\dim X-2$, which are called {\sf del Pezzo} 
and {\sf Mukai} manifolds, respectively, have been classified 
(\cite{Fub}, \cite{Mu}, \cite{Me}).\\
The method used for those cases (i.e. proving that
the linear sistem $|L|$ contains a smooth divisor and constructing a ladder
down to the known cases of lower dimensional varieties)
does not work for Fano manifolds of index $\dim X-3$, since there
are no results on the existence of a (smooth) divisor in the linear
system $|L|$ and, most of all, the classification of Fano fourfolds
is very far from being known.\\
Nevertheless some classification results for Fano manifold of index $\dim X-3$ 
and Picard number greater than one are known: by the classification
of Fano manifolds of middle index and Picard number greater than one
obtained by Wi\'sniewski and other authors (see \cite{Wimi} for a survey on these results)
we have the complete classification of Fano manifolds of index $\dim X-3$, Picard number
greater than one and dimension greater than or equal to six.\\
Roughly speaking, apart from $\pd \times \pd \times \pd$, these
varieties have Picard number two, and thus two extremal elementary contractions,
and the classification is obtained by a careful study of these
contractions and their interplay.\\
Actually, by a theorem of Wi\'sniewski \cite{Wimu}, there are no Fano manifolds
of index $\dim X-3$ and dimension greater than eight; this theorem
is a particular case of a conjecture of Mukai  
relating the pseudoindex, the dimension and the Picard number of a Fano manifolds:
$$\rho_X(i_X -1) \le \dim X.$$
In \cite{ACO} it was proved that the conjecture 
holds for Fano manifolds of dimension five (for lower dimensional cases the
result was already known).\\
However, the information on the Picard number when $\rho_X \ge 3$ is not
enough to decide the number and type of the extremal contractions of the variety,
i.e. to understand the structure of the cone of curves $\cone(X)$, result
that was achieved for Fano fivefolds of pseudoindex greater than one in \cite{CO}.\\
The present paper is intended as a first step in going from the table of the cones
given in \cite{CO} to the actual classification of Fano fivefolds of
index two, and it deals with ruled Fano fivefolds, i.e. with
triples $(X,Y, \E)$ constituted by a Fano fivefold $X$ of index two,
a smooth variety $Y$ of dimension four and a rank two vector bundle $\E$ over 
$Y$ such that $X= \proj_Y(\E)$.\par
\smallskip
The paper is organized as follows: in section \ref{back} we collect 
basic material concerning Fano-Mori contractions, 
families of rational curves and Fano manifolds; section \ref{ruled} is dedicated
to $\proj^{r-1}$-ruled Fano manifolds of index $r$, i.e. triples as above
where $\rk \E= r_X=r$, relating the extremal contractions of $X$ and $Y$.\\
Section \ref{recprod} contains some criteria to establish if a $\proj^{r-1}$-ruled
Fano manifold of index $r$ is a product of another Fano manifold
of index $r$ with a projective space $\proj^{r-1}$.\\ 
In section \ref{mafi} we begin with the classification problem; 
as already showed by the table of the cones in \cite{CO}, the greater is 
the Picard number, the easier the classification becomes; this allows us 
to treat the cases $\rho_X \ge 4$ in a broader context, proving two 
general results on Fano manifolds with large Picard number and only 
(or almost only) fiber type contractions (propositions \ref{fibertype} and \ref{mixedtype}).\\
The following two sections are dedicated to the case $\rho_X=3$, and we prove the
following

\begin{theorem}\label{rho3} Let $(X,Y, \E)$ be a ruled Fano fivefold of index two 
with $\rho_X \ge 3$; then either $X$ is a product $\pu \times Y$, with
$Y$ a Fano fourfold of index two and $\rho_Y=2$ (for a classification of these 
manifolds see \cite{WiF4}) or $X$ is one of the following:
\begin{enumerate}
\item  $X \simeq Bl_p(\proj^4) \times_\pt Bl_p(\proj^4)$;
\item $X \simeq Bl_S(Bl_p(\proj^5))$ with $S$ the strict trasform of a plane $\ni p$;
\item the blow up of $\proj^5$ in two non meeting planes;
\item the blow up of a cone in $\proj^9$ over the Segre embedding $\pd \times \pd \subset \proj^8$
along its vertex;
\item the blow up of a general member of $\Ol(1,1) \subset \pd \times \proj^4$ along
a two dimensional fiber of the second projection.
\end{enumerate}
In these cases the corresponding pairs $(Y,\E)$ are, respectively,
\begin{enumerate}
\item $(Bl_p(\proj^4), 2H +E \oplus 3H +E)$, $E$ exceptional divisor and
 $H$ pullback on $Y$ of $\Ol_\pt(1)$;
\item $(Bl_l(\proj^4), 2H-E \oplus 3H-E)$, $E$ exceptional divisor and
 $H$ pullback on $Y$ of $\Ol_{\proj^4}(1)$;
\item $(\pd \times \pd, \Ol(1,2) \oplus \Ol(2,1))$;
\item $(\pd \times \pd, \Ol(1,1) \oplus \Ol(2,2))$;
\item $(\proj_\pd(T\pd(-1)\oplus\Ol_\pd)\subset\pd \times \pt, \Ol(1,1) \oplus \Ol(1,2))$.
\end{enumerate}
\end{theorem}

The last section contains the case $\rho_X=2$, in which we have the following

\begin{theorem}\label{rho2} Let $(X,Y, \E)$ be a ruled Fano fivefold of index two 
with $\rho_X = 2$; then either $X$ is a product $\pu \times \mathbb Q^4$, or
$\pu \times Y$ with $Y$ a Mukai fourfold of Picard number one (see \cite{Mu})
or $X$ is one of the following:
\begin{enumerate}
\item $\proj_{\proj^4}(\Ol_{\proj^4} \oplus \Ol_{\proj^4}(a))$, with $a=1$ or $a=3$;
\item $\proj_{\mathbb Q^4}(\Ol_{\mathbb Q^4} \oplus \Ol_{\mathbb Q^4}(2))$;
\item $\proj_{V_d}(\Ol_{V_d} \oplus \Ol_{V_d}(1))$, with $V_d$ a del Pezzo fourfold of degree 
$d=1,\dots, 5$;
\item a general divisor in the linear system $|2 \xi|$ in
$\proj_\pt(\Omega\pt(3) \oplus \Ol(1))$;
\item in $G(1,4) \times \proj^4$, the intersection of two divisors in the linear 
system $|\Ol(1,0)|$ with the flag variety of point and lines in $\proj^4$;
\item a $\pu$-bundle over a Fano fourfold of index one and pseudoindex two or three.
\end{enumerate}
\end{theorem}

Our classification is effective, apart from case (6) of theorem \ref{rho2};
we point out that it is not known whether a Fano fourfold as in case (6) (i.e. a Fano fourfold
of Picard number one without a line) exists or not, and its existence (or non existence)
constitutes a very hard problem.


\section{Background material}\label{back}

\subsection{Extremal contractions} 
Let $X$ be a smooth complex Fano variety
of dimension $n$ and let $K_X$ be its canonical divisor.
By Mori's {\sf Cone Theorem} the cone of effective 1-cycles, which is contained in
the $\mathbb R$-vector space of 1-cyles modulo numerical equivalence,
$\cone(X) \subset N_1(X)$, is polyhedral; a face of $\cone(X)$ is
called an {\sf extremal face} and an extremal face of dimension one
is called an {\sf extremal ray}.\\
To an extremal face $\sigma \subset \cone(X)$ is associated a morphism with connected fibers 
$\f_\sigma:X \to Z$ onto a normal variety, morphism which contracts the curves 
whose numerical class is in $\sigma$; $\f_\sigma$ is called an {\sf extremal contraction}
or a {\sf Fano-Mori contraction}, while a Cartier divisor $H$ such that 
$H = \f_\sigma^*A$ for an ample divisor $A$ on $Z$
is called a {\sf supporting divisor} of the map $\f_\sigma$ (or of the face $\sigma$).\\
An extremal contraction associated to an extremal ray is called an 
{\sf elementary contraction};
an extremal ray $R$ is called {\sf numerically effective}, and the
associated contraction is said to be of {\sf fiber type}, if $\dim Z < \dim X$;
otherwise the ray is called {\sf non nef} and the contraction is {\sf birational};
the terminology is due to the fact that, if $R$ is a non nef ray, there 
exists an irreducible divisor which has negative intersection number with curves in $R$.\\
We usually denote with $\Exc(\f_\sigma):= \{x \in X\ |\ \dim \f_\sigma^{-1}(\f_\sigma(x)) > 0\}$
the {\sf exceptional locus} of $\f_\sigma$; if $\f_\sigma$ is of fiber type then, 
of course, $\Exc(\f_\sigma)=X$.\\
If the codimension of the exceptional locus of an elementary birational contraction
is equal to one, the ray and the contraction are called {\sf divisorial}, otherwise they are 
called {\sf small}.

\begin{definition} An elementary fiber type extremal contraction 
$\f: X \to Z$ is called a {\sf scroll} (respectively a {\sf quadric fibration}) 
if there exists a $\f$-ample line bundle $L \in \pic(X)$ such that
$K_X+(\dim X-\dim Z+1)L$ (respectively $K_X+(\dim X-\dim Z)L$)
is a supporting divisor of $\f$; we will call {\sf conic fibration}
a quadric fibration such that $\dim X-\dim Z=1$.\\
An elementary fiber type extremal contraction $\f:X \to Z$ onto a smooth
variety $Z$ is called a $\mathbb P${\sf-bundle} (respectively
{\sf quadric bundle}) if there exists a vector bundle $\E$ of rank 
$\dim X-\dim Z+1$ (respectively of rank $\dim X-\dim Z+2$) on $Z$ such that 
$X \simeq \proj(\E)$ (respectively there exists an embedding of $X$ over $Z$ as a divisor of
$\proj(\E)$ of relative degree 2); we will call {\sf conic bundle}
a quadric bundle such that $\dim X-\dim Z=1$.\\
An equidimensional scroll is a projective bundle by \cite[Lemma 2.12]{Fu1},
while an equidimensional quadric fibration is a quadric bundle by
\cite[Theorem B]{ABW3}.\\ 
Some special scroll contractions arise from projectivization of
B\v anic\v a sheaves (cfr. \cite{BW}); in parti\-cu\-lar,
if $\f:X \to Z$ is a scroll such that every fiber has dimension $\le \dim X- \dim Z+1$,
then $Z$ is smooth and $X$ is the projectivization of a B\v anic\v a sheaf on $Z$
(cfr. \cite[Proposition 2.5]{BW}); we will call these contractions
{\sf special B\v anic\v a scrolls}.
\end{definition}

\subsection{Families of rational curves} 
For this subsection our main reference is \cite{Kob}, 
with which our notation is coherent.
Let $X$ be a normal projective variety and let $\om(\pu,X)$ be the scheme parametrizing
morphisms $f: \pu \to X$; let $\om_{bir}(\pu,X) \subset \om(\pu,X)$ be the open subscheme 
corresponding to those morphisms which are birational onto their image, and let
$\om^n_{bir}(\pu,X)$ be its normalization; the group $\Aut(\pu)$ acts on 
$\om^n_{bir}(\pu,X)$ and the quotient exists.\par

\smallskip
\begin{definition}
The space $\ratcurves$ is the quotient of $\om^n_{bir}(\pu,X)$ by $\Aut(\pu)$, and the space
$\Univ(X)$ is the quotient of the product action of $\Aut(\pu)$ on
$\om^n_{bir}(\pu, X)\times\pu$.
\end{definition}

\begin{definition} \label{Rf}
A {\sf family of rational curves} is an irreducible component
$V \subset \ratcurves$.\\
Given a rational curve $f:\pu \to X$, we will call a {\sf family of
deformations} of $f$ any irreducible component $V \subset
\ratcurves$ containing the equivalence class of $f$.\\
Given a family $V$ of rational curves, we have the following basic diagram
$$
\xymatrix{p^{-1}(V)=:U  \ar[r]^(.65){i} \ar[d]^{p} & X\\
 V & & }
$$
where $i$ is the map induced by the evaluation $ev:\om^n_{bir} (\pu, X)\times \pu \to X$
and $p$ is the $\pu$-bundle induced by the projection $\om^n_{bir}(\pu,X) \times \pu \to
\om^n_{bir}(\pu,X)$.
We define $\loc(V)$ to be the image of $U$ in $X$;
we say that $V$ is a {\sf covering family} if $\overline{\loc(V)}=X$.\\
If $L \in \pic(X)$ is a line bundle, we will denote by $L \cdot V$ the intersection
number of $L$ and a general member of the family $V$.
Finally, given a family $V \subseteq \ratcurves$, we denote by $V_x$ the subscheme
of $V$ parametrizing rational  curves passing through $x$.
\end{definition}

\begin{definition} 
Let $V$ be a family of rational curves on $X$. Then $V$ is {\sf unsplit} if it is proper. 
\end{definition}

\begin{example}\label{unex} 
Let $R_i$ be an extremal ray and $C_i$ a curve whose numerical class belongs to $R_i$ 
and whose anticanonical degree is minimal among curves whose class is in $R_i$; 
$C_i$ is often called a {\sf minimal extremal rational curve}.\\
Denote by $R^i$ an irreducible component of $\ratcurves$ containing $C_i$; then
the family $R^i$ is unsplit: indeed, if $C_i$ degenerates into a reducible cycle, 
its components must belong to the ray $R_i$, since $R_i$ is extremal; but in $R_i$
the curve $C_i$ has the minimal intersection with the anticanonical
bundle, hence this is impossible. 
\end{example}

\begin{proposition}\cite[IV.2.6]{Kob}\label{iowifam} 
Let $X$ be a smooth projective variety
and $V$ a family of rational curves.
Assume 
that $V$ is unsplit and $x$ is any point in $\loc(V)$. Then
 \begin{itemize}
      \item[(a)] $\dim X -K_X \cdot V \le \dim \loc(V)+\dim \loc(V_x) +1$;
      \item[(b)] $-K_X \cdot V \le \dim \loc(V_x)+1$.
   \end{itemize}
\end{proposition}

This last proposition, in case $V$ is the unsplit family of deformations of a minimal extremal
rational curve, gives the {\sf fiber locus inequality}:

\begin{proposition}\label{fiberlocus} Let $\f$ be a Fano-Mori contraction
of $X$ and let $E = \Exc(\f)$ be its exceptional locus;
let $S$ be an irreducible component of a (non trivial) fiber of $\f$. Then
$$\dim E + \dim S \geq \dim X + l -1,$$
where 
$$l =  \min \{ -K_X \cdot C\ |\  C \textrm{~is a rational curve in~} S\}.$$
If $\f$ is the contraction of a ray $R$, then $l(R):=l$ is called the {\sf length of the ray}.
\end{proposition}

\noindent Let $X$ be a smooth variety, $V^1, \dots, V^k$ unsplit families of rational curves
on $X$ and $Z \subset X$.\par

\begin{definition}
We denote by $\loc(V^1, \dots, V^k)_Z$ the set of points that can be joined to $Z$ by 
a connected chain of $k$ cycles belonging \underline{respectively} to the families 
$V^1, \dots, V^k$.\\
We denote by $\cloc_m(V^1, \dots, V^k)_Z$ the set of points that can be joined to $Z$ by a
connected chain  of at most $m$ cycles belonging to the families $V^1, \dots, V^k$.
\end{definition}

\begin{definition}
We define a relation of {\sf rational connectedness with respect to $V^1, \dots, V^k$}
on $X$ in the following way: $x$ and $y$ are in rc$(V^1,\dots,V^k)$-relation if there
exists a chain of rational curves in $V^1, \dots ,V^k$ which joins $x$ and $y$, i.e.
if $y \in \cloc_m(V^1, \dots, V^k)_x$ for some $m$.
\end{definition}

To the rc$(V^1,\dots,V^k)$-relation we can associate a fibration, 
at least on an open subset.

\begin{theorem}\cite{Cam},\cite[IV.4.16]{Kob} \label{rcvfibration}
There exist an open subvariety $X^0 \subset X$ and a proper morphism with connected fibers
$\pi:X^0 \to T^0$ such that
   \begin{itemize}
      \item[(a)] the rc$(V^1,\dots,V^k)$-relation restricts to an equivalence relation 
	  on $X^0$;
      \item[(b)] the fibers of $\pi$ are equivalence classes for the 
	  rc$(V^1,\dots,V^k)$-relation;
      \item[(c)] for every $t \in T^0$ any two points in $\pi^{-1}(t)$ can be connected 
	  by a chain of at most \linebreak $2^{\dim X - \dim T^0}-1$ cycles in $V^1, \dots, V^k$.
   \end{itemize}
\end{theorem}

\begin{definition}
In the above assumptions, if $\pi$ is the constant map, we will say that $X$ is 
{\sf rc$(V^1,\dots,V^k)$-connected}.
\end{definition}

For other properties of $\loc(V^1, \dots, V^k)_Z$ and $\cloc_m(V^1, \dots, V^k)_Z$
we refer to \cite{ACO} and \cite{CO}.

\subsection{Fano manifolds and projective bundles}

\begin{lemma}\label{bundlebase}
Let $X$ be a Fano manifold and $p:X \to Y$ an elementary contraction
onto a smooth variety such that every fiber of $p$ is a projective
space of dimension $r$. Denote by $R_\E$ the extremal ray of $\cone(X)$ 
corresponding to $p$. Then
\begin{itemize}
\item[(a)] $Y$ is a Fano manifold with pseudoindex $i_Y \ge i_X$;
\item[(b)] if $i_Y=i_X$ and $f:\pu \to Y$ is a rational curve of degree $i_Y$, then 
$f^*\E \simeq \Ol_\pu(a)^{\oplus r+1}$;
\item[(c)] if $\cone(X)=\langle R_\E, R_1, \ldots, R_k \rangle$, then 
$\cone(Y)=\langle p(R_1), \ldots, p(R_k) \rangle$.
\end{itemize} 
\end{lemma}

{\bf Proof.} \quad $Y$ is a Fano manifold by \cite[Corollary 2.9]{KoMiMo}; 
the assertion on the pseudoindex and part (b) are proved in \cite[Lemma 2.5]{BCDD}, 
while part (c) is contained in the proof of \cite[Lemma 3.1]{Wi1}.\qed \par

\begin{lemma}\label{pubundle}
Let $X$ be a Fano manifold of pseudoindex $i_X \ge 2$ and let $\f:X \to Y$ be an elementary 
contraction which is equidimensional with one dimensional fibers.
Then there exists a rank two vector bundle 
$\E$ on $Y$ such that $X=\proj_Y(\E)$.
\end{lemma}

{\bf Proof.} \quad By \cite[Theorem 3.1 (ii)]{An} $Y$ is smooth and $\f:X \to Y$ is
a conic bundle. It follows that $-K_X \cdot f=2$ for every fiber $f$ of $\f$, 
therefore $f$ can not be reducible or nonreduced, being $i_X \ge 2$.
By lemma \ref{bundlebase} (a) $Y$ is a Fano manifold; in particular its Brauer
group is trivial, hence there exists a rank two vector bundle 
$\E$ on $Y$ such that $X=\proj_Y(\E)$.\qed \par
\medskip
The fact that cone of curves of a Fano manifold is polyhedral and generated by a finite
number of extremal rays easily leads to the following

\begin{lemma}\label{diveff} \cite[Lemme 2.1]{BCW} Let $X$ be a Fano manifold and $D$ an 
effective divisor on $X$. Then there exists an extremal ray $R \subset \cone(X)$
such that $D \cdot R >0$.
\end{lemma}

which, combined with lemma \ref{pubundle}, gives

\begin{corollary} \label{onedim} Let $X$ be a Fano manifold of pseudoindex $i_X \ge 2$,
$R \subset \cone(X)$ an extremal ray and $D$ an effective divisor on $X$ 
such that no curve in $D$ has numerical class belonging to $R$.
If $D \cdot R>0$, then the contraction associated to $R$, $\f_R:X \to Y$ is a $\pu$-bundle.
\end{corollary}

{\bf Proof.} \quad Let $F$ be any fiber of $\f_R$; the intersection $D \cap F$
has to be zero dimensional, otherwise $D$ would contain a curve whose numerical
class is in $R$. It follows that $\f_R$ is equidimensional with one dimensional
fibers and we can apply lemma \ref{pubundle}. \qed \par
\medskip
The following lemma will be of frequent use in our proofs:

\begin{lemma}\label{fano3fold}
Let $T$ be a smooth threefold of Picard number one, $\F$ a rank two vector bundle 
on $T$ and $Y=\proj_T(\F)$; assume that $Y$ is a Fano manifold of pseudoindex $i_Y \ge 2$. 
Then, if $Y$ is not a product $Y = \pu \times T$, we have
either $T \simeq \pt$ or $T \simeq \qt$.
\end{lemma}

{\bf Proof.} \quad By lemma \ref{bundlebase} (a), $T$ is a Fano threefold of pseudoindex
$i_T \ge i_Y \ge 2$; in particular, by the classification of Fano threefolds, 
$T$ admits an unsplit covering family $V_T$ of rational curves of degree $i_T$.\\
If $i_Y=i_T$, then, by lemma \ref{bundlebase} (b),
the restriction of $\F$ to any curve of $V_T$ splits as $\Ol_\pu(a) \oplus \Ol_\pu(a)$,
so, by \cite[Proposition 1.2]{AWI}, $\F$ is decomposable and $Y \simeq \pu \times T$.\\
Otherwise $i_T \ge 3$ and, by the classification of Fano threefolds, 
either $T \simeq \qt$, or $T \simeq \pt$.\qed\par
\medskip
Finally we prove two lemmata which ensure that, in some cases, a fibration
in projective spaces is a projective bundle.

\begin{lemma}\label{curve} 
Let $p:Y \to B$ be a morphism from a smooth variety 
to a smooth curve, such that $\rho(Y/B)=1$ and the general fiber of $p$ is a projective space; 
then there exists a vector bundle ${\mathcal F}$ of rank = $\dim Y$ on $B$ such that 
$Y=\proj_B({\mathcal F})$ and $p$ is the natural projection. 
\end{lemma}

{\bf Proof.} \quad Over an open Zariski subset $U$ of $B$ the morphism $p$
is a projective bundle; indeed over a curve $C$ a fibration in projective spaces
is a projective bundle, since the obstruction lies in $H^2(C,\Ol^*)=0$ (see \cite{El}).
By taking the closure in $Y$ of a
hyperplane section of $p$ defined over the open
set $U$ we get a global relative hyperplane section divisor (we use
$\rho(Y/B)=1$) hence $p$ is a projective bundle globally by \cite[Lemma 2.12]{Fu1}.\qed

\begin{lemma}\label{brauer}
Let $X$ be a Fano manifold  and $p:X \to S$ be an elementary
contraction associated to an extremal ray of length $\dim X-1$ onto a surface $S$. 
Then $S$ is smooth and there exists a rank $\dim X-1$ vector bundle $\F$ over $S$
such that $X=\proj_{S}(\F)$.
\end{lemma}

{\bf Proof.} Since $p$ is elementary and $\dim S=2$ then $p$ is equidimensional;
by \cite[Corollary 1.4]{AW} $S$ is smooth.\\
By adjunction the general fiber of $p$ is a projective space of dimension $\dim X-2$;
over a general hyperplane section of $S$, $\f$ is a  projective bundle by 
lemma \ref{curve}, whence the locus over which the fiber is not a projective space
is discrete in $S$.
We can apply \cite[Lemma 3.3]{AnMe} and \cite[Lemma 2.12]{Fu1} to obtain that 
every fiber of $\f$ is a projective space.
The surface $S$ is dominated by a Fano manifold, hence is rationally connected;
therefore $H^2(S,\Ol^*)=0$ and the Brauer group of $S$ is trivial.
This implies the existence of a rank $\dim X-1$ vector bundle $\F$ over
$S$ such that $X=\proj_S(\F)$.\qed  


\section{$\proj^{r-1}$-ruled Fano manifolds: general properties}\label{ruled}

\begin{definition}\label{rfm}
Let $Y$ be a smooth variety of dimension $n$, let $\E$ be a vector bundle
of rank $r$ on $Y$ and let $X= \proj_Y(\E)$ be the projectivization of $\E$;
assume moreover that $X$ is a Fano manifold.
We will call a triple $(X,Y, \E)$ as above a $\proj^{r-1}$-{\sf ruled Fano manifold}; 
if $r=2$, we will call for short $(X,Y,\E)$ a 
{\sf ruled Fano manifold}.
\end{definition}

\begin{definition}\label{rfmi}
Let $(X,Y,\E)$ be a $\proj^{r-1}$-ruled Fano manifold verifying one
of the following
\begin{enumerate}
\item[1)] $X$ has index $r$;
\item[2)] $K_Y + \det \E' =\Ol_Y$, with $\E'$ an ample twist of $\E$.
\end{enumerate}
We will call such a triple a $\proj^{r-1}$-{\sf ruled Fano manifold 
of index $r$}; if $r=2$, we will call for short $(X,Y,\E)$ a 
{\sf ruled Fano manifold of index two}.\\
From now on, unless otherwise stated, we 
will assume that $\proj^{r-1}$-{\sf ruled Fano manifold of index $r$} $(X,Y,\E)$ are
normalized, i.e. $\E$ is ample and $K_Y + \det \E =\Ol_Y$.
\end{definition}

\begin{remark} The assumptions 1) and 2) are equivalent.
\end{remark}

{\bf Proof.} \quad Let us show first that 1) $\Rightarrow$ 2); let $H \in \pic(X)$ be
the (unique) line bundle such that $-K_X=rH$; by adjunction, if $l$ is a line 
in a fiber of the projection $p:X \to Y$, then $r= -K_X \cdot l =rH \cdot l$, so 
$H$ restricts to $\Ol_{\proj^{r-1}}(1)$ on the fibers of $p$.
Therefore $p_*H$ is an ample vector bundle of rank $r$, $\E'$, which differs
from $\E$ by a twist with a line bundle in $\pic(Y)$ and, by the canonical bundle
formula 
$$\Ol_X=K_X+rH=p^*(K_Y+ \det \E'),$$
hence $K_Y+ \det \E'=\Ol_Y$.\\
Assume now that 2) holds; for a suitable ample twist $\E' = \E \otimes L$,
we have $K_Y + \det \E' =\Ol_Y$, therefore, by the canonical bundle formula,
$$K_X +r\xi_{\E'} = p^*(K_Y+\det \E')= \Ol_X,$$
whence $-K_X=r\xi_{\E'}$ and $X$ is a Fano manifold of index $r$. \qed \par

\begin{proposition} \label{corray} Let $(X,Y,\E)$ be a $\proj^{r-1}$-ruled
Fano manifold and denote by $R_\E$ the extremal ray in $\cone(X)$ 
associated to the bundle projection $p:X \to Y$. There is a one-to-one correspondence
$$
\xymatrix@1{\txt{$\Bigg \{$} \txt{Extremal rays of $\cone(X)$ spanning\\
a two dimensional face with $R_\E$} \txt{$\Bigg \}$}}
{\xymatrix@1{
\ar@/^0.5pc/[rr]<1ex>^-{\alpha_X}
&&\ar@/^0.5pc/[ll]<1ex>^-{\alpha_Y}
}}
\xymatrix@1{\txt{$\bigg \{$} \txt{Extremal rays of $\cone(Y)$} \txt{$\bigg \}$}}.
$$
If $\theta \subset \cone(Y)$ and $\vartheta \subset \cone(X)$ are corresponding
rays, then we will call them {\sf fellow} rays.
\end{proposition}

{\bf Proof.} \quad
Let $\theta$ be an extremal ray of $\cone(Y)$ and denote by $\f_\theta:Y \to W$ the associated
elementary contraction; then $\rho(X/W)=2$ and $-K_X$ is $(\f_\theta \circ p)$-ample, so
$\f_\theta \circ p: X \to W$ is the contraction of a two dimensional extremal face
$\sigma \subset \cone(X)$ containing $R_\E$. Let $\vartheta$ be the
extremal ray in $\sigma$ different from $R_\E$; we set $\alpha_Y(\theta)=\vartheta$.\par
\medskip
On the other hand, if $\vartheta$ is an extremal ray of $\cone(X)$ such that 
$\sigma= \langle R_\E, \vartheta \rangle$ is an extremal face, then 
the contraction $\psi_\sigma:X \to W$ factors 
both through the contraction $p$ of $R_\E$ and through the contraction $\psi_\vartheta:X \to Z$
of $\vartheta$, hence we have a commutative diagram
\begin{equation} \label{gc}
\xymatrix@=40pt{X \ar[rd]^{\psi_\sigma} \ar[r]^{\pst} \ar[d]_{p} & Z \ar[d]^{p'}\\ 
Y \ar[r]_\ft  &W }
\end{equation}

Since $Y$ is a Fano manifold and $\f_\theta$ is a surjective 
morphism with connected fibers, we have that $\f_\theta$ is an extremal contraction; 
moreover, being $\rho(Y/W)=1$, the contraction is elementary, thus it corresponds 
to an extremal ray $\theta$. Setting $\alpha_X(\vartheta)=\theta$ we have the
desired bijection.\qed \par

\setcounter{equation}{0}
\begin{lemma}\label{fiberdim}
Let $(X,Y,\E)$ be a $\proj^{r-1}$-ruled Fano manifold and let 
$\theta \subset \cone(Y)$ and $\vartheta \subset \cone(X)$ 
be two fellow rays with associated extremal contractions $\f_\theta:Y \to W$ and
$\psi_\vartheta:X \to Z$, with exceptional loci $\Exc(\f_\theta)$
and $\Exc(\psi_\vartheta)$ respectively. Then
\begin{equation}\label{fiberdim1}
p(\Exc(\psi_\vartheta)) \subset \Exc(\f_\theta).
\end{equation}
Moreover, if $x$ is a point in $\Exc(\psi_\vartheta)$, $(F_\psi)_x$
is the fiber of $\psi_\vartheta$ through $x$ and $(F_\f)_{p(x)}$ is the 
fiber of $\f_\theta$ through $p(x)$, we have 
\begin{equation}\label{fiberdim2}
\dim (F_\psi)_x =\dim p((F_\psi)_x)\le \dim (F_\f)_{p(x)}.
\end{equation}
Finally, if $x_1$ is a point in $p^{-1}(p(x)) \cap \Exc(\psi_\vartheta)$
and $(F_\psi)_{x_1}$ is the fiber of $\pst$ through $x_1$, then
\begin{equation}\label{fiberdim3}
p(F_\psi)_{x_1} \subset (F_\f)_{p(x)}.
\end{equation}
\end{lemma}

{\bf Proof.} \quad The statements follows from the commutativity of diagram
\ref{gc} and the fact that the projection $p$, being the contraction of an extremal
ray different from $\vartheta$, is finite to one on the fibers of $\pst$. \qed \par

\begin{corollary} \label{excloc} Under the assumptions of lemma \ref{fiberdim},
if $\psi_\vartheta$ is of fiber type then also $\f_\theta$ is of fiber type, 
while if $\f_\theta$ is birational then also $\psi_\vartheta$ is birational.
\end{corollary}

\begin{lemma}\label{supporting} Let $(X,Y,\E)$ be a $\proj^{r-1}$-ruled
Fano manifold of index $r$ and let $\theta \subset \cone(Y)$ and $\vartheta \subset \cone(X)$ 
be two fellow rays with associated extremal contractions $\f_\theta:Y \to W$ and
$\psi_\vartheta:X \to Z$. Then there exist an ample vector bundle $\E_\Theta$
on $Y$ and an ample line bundle $L \in \pic(X)$ such that $\f_\theta$ is supported by 
$K_Y + \det \E_\Theta$ and $\psi_\vartheta$ is supported by $K_X + rL$.
\end{lemma}

{\bf Proof.} \quad  Pick two ample line bundles
$A \in \pic(W)$ and $B \in \pic(Z)$. Set $\E_\Theta=\E \otimes \f_\theta^*A$;
we have $K_Y +\det \E_\Theta=r\f_\theta^*A$, so we have only to prove the ampleness
of $\E_\Theta$.\\
The tautological line bundle associated to $\E_\Theta$ on
$\proj(\E_\Theta)=\proj(\E)=X$ is 
$$\xi_\Theta= \xi_\E + p^*(\f_\theta^*A),$$ 
hence it is ample, being the sum of an ample line bundle and a nef one.\par
\medskip
To prove the second statement observe that $K_X +r\xi_{\Theta} 
=p^*(K_Y +\det\E_\Theta) = r(p^*(\f_\theta^*A))$; therefore, if 
$L:=\xi_\Theta +\psi_\vartheta^*B$, we have
$$K_X +r L= r(p^*(\f_\theta^*A)+\psi_\vartheta^*B)=
r\psi_\vartheta^*(p'^*A +B).$$ 
Moreover $L$ is ample, being the sum
of an ample line bundle and a nef one.\qed \par
\medskip
We now analyze some cases in which $\f_\theta$ is a special
contraction (projective bundle, smooth blow-up, special
\ba scroll), describing the structure of the corresponding contraction $\pst$.

\begin{proposition}\label{spcontr} Let $(X,Y,\E)$ be a $\proj^{r-1}$-ruled
Fano manifold of index $r$; let $\theta \subset \cone(Y)$ and $\vartheta \subset \cone(X)$
be two fellow rays and let $\f_\theta:Y \to W$ and $\psi_\vartheta:X \to Z$ 
be the associated contractions. \nolinebreak Then
\begin{enumerate}
\item[(a)] if $\f_\theta$ is a $\proj^{r-1}$-bundle, then $\psi_\vartheta$ 
is a $\proj^{r-1}$-bundle;
\item[(b)] if $\f_\theta$ is the blow up of a smooth subvariety of $W$ of codimension
$r+1$, then $\psi_\vartheta$ is the blow up of a smooth subvariety of $Z$ of codimension
$r+1$.
\end{enumerate} 
In both cases, if $H \in \pic(Y)$ is a line bundle which restricts to $\Ol_{\proj}(1)$
on the fibers of $\f_\theta$, then $\E \otimes H^{-1}=\f_\theta^* \E'$, 
where $\E'$ is a rank $r$ vector bundle on $W$, and $Z=\proj_W(\E')$.
\end{proposition}

{\bf Proof.} \quad Denote by $l(\vartheta)$ the length of the extremal ray $\vartheta$;
since $X$ is a Fano manifold of index $r$ we have $l(\vartheta) \ge r$.\\
In case (a), if $x \in X$ is any point in $\Exc(\vartheta)$, $(F_\psi)_x$ 
is the fiber of $\psi_\vartheta$ through $x$ and
$(F_\f)_{p(x)}$ is the fiber of $\f_\theta$ through $p(x)$, by proposition
\ref{fiberlocus} and formula \ref{fiberdim2} we have 
$$r-1 \le l(\vartheta)-1 \le \dim (F_\psi)_x \le \dim (F_\f)_{p(x)} =r-1,$$
so $\psi_\vartheta$ is an equidimensional contraction with $(r-1)$-dimensional
fibers (and thereby of fiber type, by proposition \ref{fiberlocus}).
By lemma \ref{supporting}, there exists an ample $L \in \pic(X)$ such that 
$\psi_\vartheta$ is supported by $K_X + rL$, and we conclude by \cite[Lemma 2.12]{Fu1}.\\
In case (b), by corollary \ref{excloc}, since $\f_\theta$ is birational, also $\psi_\vartheta$ 
is birational. 
Then, if $x \in X$ is any point in $\Exc(\vartheta)$, $(F_\psi)_x$ 
is the fiber of $\psi_\vartheta$ through $x$ and
$(F_\f)_{p(x)}$ is the fiber of $\f_\theta$ through $p(x)$, by proposition
\ref{fiberlocus} and formula \ref{fiberdim2} we have 
$$r \le l(\vartheta) \le \dim (F_\psi)_x \le \dim (F_\f)_{p(x)} = r,$$
thus $\psi_\vartheta$ is equidimensional
with fibers of dimension $r$ and, by lemma \ref{supporting}, it is 
supported by $K_X+rL$, for some ample $L \in \pic(X)$; therefore we can apply
\cite[Theorem 4.1]{AWD} to conclude.\\
In both cases the extremal ray $\theta$ has length $r$, hence $r \ge i_Y$;
by lemma \ref{bundlebase} (a) we have $i_Y \ge i_X$ and, recalling
that the pseudoindex $i_X$ is greater or equal than the index $r_X=r$, we have $i_X \ge r$.
We conclude that $i_Y = i_X = r$.\\ 
By lemma \ref{bundlebase} (b), for every line $l$ 
in every fiber of $\f_\theta$ we have $\E_l \simeq \Ol_\pu(1)^{\oplus r}$, hence,
if $H \in \pic(Y)$ is a line bundle which restricts to $\Ol_{\proj}(1)$ on the fibers
of $\f_\theta$, the vector bundle $\E \otimes H^{-1}$ is trivial on every fiber, so
it is the pullback of a rank $r$ vector bundle $\E'$ on $W$.
It is now easy to prove that the induced map $\proj_Y(\f_\theta^*\E')= X \to \proj_W(\E')$
is just $\psi_\vartheta$, whence $Z= \proj_W(\E')$. \qed

\begin{proposition}\label{spcontr2} Let $(X,Y,\E)$ be a $\proj^{r-1}$-ruled
Fano manifold of index $r$; let $\theta \subset \cone(Y)$ and $\vartheta \subset \cone(X)$ 
be two fellow rays and let $\f_\theta:Y \to W$, $\psi_\vartheta:X \to Z$ 
be the associated contractions. \nolinebreak Then
\begin{enumerate}
\item[(a)] if $\f_\theta$ is a $\proj^{r}$-bundle and $\psi_\vartheta$ is of fiber type, 
then $\psi_\vartheta$ is a $\proj^{r-1}$-bundle;
\item[(b)] if $\f_\theta$ is a $\proj^{r}$-bundle and $\psi_\vartheta$ is birational, 
then $\psi_\vartheta$ is the blow up of a codimension $r+1$ subvariety of $Z$.
\end{enumerate} 
Moreover, in case (a), if $H \in \pic(Y)$ is a line bundle which restricts to 
$\Ol_{\proj^r}(1)$ on the fibers of $\ft$, then $p^*H$ restricts to $\Ol_{\proj^{r-1}}(1)$
on the fibers of $\pst$; in case (b), the divisor $\Exc(\pst)$ restricts to
$\Ol_{\proj^{r-1}}(1)$ on the fibers of $p$.
\end{proposition}
\nopagebreak[4]
{\bf Proof}.\quad Let $\sigma= \langle R_\E, \vartheta \rangle \subset \cone(X)$ and let
$\psi_\sigma:X \to W$ be the contraction associated to the face $\sigma$, which can be 
factored both as $\ft \circ p$ and as $p' \circ \pst$:

$$
\xymatrix@=40pt{X \ar[rd]^{\psi_\sigma}  \ar[r]^{\pst} \ar[d]_{p}  & Z \ar[d]^{p'}\\
Y \ar[r]_\ft  &W }
$$

A fiber $F_\sigma$ of $\psi_\sigma$ can thus be viewed as the inverse image 
via $p$ of a fiber $F_\theta\simeq \proj^r$ of $\ft$, 
$F_\sigma \simeq \proj_{F_\theta}(\E_{|F_\theta})$.\\
The ampleness of the vector bundle $\E$ together with the fact that 
$$\det (\E_{|F_\theta}) = (\det \E)_{|F_\theta}=(-K_Y)_{|F_\theta}=\Ol_{\proj^r}(r+1)$$ 
yields that the splitting type of $\E$ on lines of $F_\theta$ is constantly
$\Ol_\pu(1)^{\oplus r-1} \oplus \Ol_\pu(2)$; by \cite{EHS}, either 
$\E_{|F_\theta} \simeq \Ol_{\proj^r}(1)^{\oplus r-1} \oplus \Ol_{\proj^r}(2)$, 
or $\E_{|F_\theta}\simeq T\proj^r$.\par
\medskip 
In case (a) $\pst$ is of fiber type, so also
its restriction to $F_\sigma=\pst^{-1}(\pst(F_\sigma))$ is a fiber type contraction, 
therefore $\E_{|F_\theta} \simeq T\proj^r$; it follows that $\pst$ 
is equidimensional and each of its fibers is $\proj^{r-1}$.
By lemma \ref{supporting}, there exists an ample $L \in \pic(X)$ such that 
$\psi_\vartheta$ is supported by $K_X + rL$, hence, by \cite[Lemma 2.12]{Fu1},
$\pst$ is a $\proj^{r-1}$-bundle over $Z$.\\ 
From this description it is clear that, if $H \in \pic(Y)$ is a line bundle which restricts to 
$\Ol_{\proj^r}(1)$ on the fibers of $\ft$, then $p^*H$ restricts to $\Ol_{\proj^{r-1}}(1)$
on the fibers of $\pst$.\par
\medskip
In case (b), if $x \in X$ is any point in $\Exc(\vartheta)$, $(F_\psi)_x$ 
is the fiber of $\psi_\vartheta$ through $x$ and
$(F_\f)_{p(x)}$ is the fiber of $\f_\theta$ through $p(x)$, by proposition
\ref{fiberlocus} and formula \ref{fiberdim2} we have 
$$r \le l(\vartheta) \le \dim (F_\psi)_x \le \dim (F_\f)_{p(x)} = r,$$
thus $\psi_\vartheta$ is equidimensional
with fibers of dimension $r$ and, by lemma \ref{supporting}, it is 
supported by $K_X+rL$, for some ample $L \in \pic(X)$; therefore, by
\cite[Theorem 4.1]{AWD} $\pst$ is the blow up of a codimension $r+1$
subvariety of $Z$.\\
Let $F_\vartheta$ be a fiber of $\pst$ and let $F_\sigma$ be the fiber
of $\psi_\sigma$ containing $F_\vartheta$; the restriction of
$\psi_\sigma$ to this fiber has a non trivial fiber of dimension $r$, therefore 
$\E_{|F_\theta} \simeq \Ol_{\proj^r}(1)^{\oplus r-1} \oplus \Ol_{\proj^r}(2)$.\\ 
It follows that $F_\sigma$ is the blow up of $\proj^{2r-1}$ along $\proj^{r-1}$
and $\Exc(\pst)_{|F_\sigma}$ is the exceptional divisor of this blow up,
hence it restricts to $\Ol_{\proj^{r-1}}(1)$ on the fibers of $p$.\qed

\begin{proposition}\label{spcontr3} Let $(X,Y,\E)$ be a $\proj^{r-1}$-ruled
Fano manifold of index $r$; let $\theta \subset \cone(Y)$ and $\vartheta \subset \cone(X)$ 
be two fellow rays and let $\f_\theta:Y \to W$ and $\psi_\vartheta:X \to Z$ 
be the associated contractions.
If $\f_\theta$ is a special \ba scroll with general fiber of dimension $r-1$, then 
also $\psi_\vartheta$ is a special \ba scroll with general fiber of dimension $r-1$.
Moreover, if $J$ is a jumping fiber of $\ft$ (i.e. a fiber of dimension $r$), 
then there is an isomorphism $f: \proj^{r-1} \times J \to  p^{-1}(J)$ 
and, for every $x \in \proj^{r-1}$,
$f(\{x \times J\})$ is a jumping fiber of $\pst$.
\end{proposition}

{\bf Proof}. \quad The general fiber of $\ft$ is $r-1$ 
dimensional, and every fiber of $\ft$ has dimension $\le r$;
using formula \ref{fiberdim2}, as in the proof of proposition \ref{spcontr}
we find that the same is true for $\pst$.\\
By lemma \ref{supporting}, the contraction $\pst$ is
supported by $K_X+rL$ for some ample $L \in \pic(X)$; we can thereby apply
\cite[Proposition 2.5]{BW} to conclude that $\pst$ is a special \ba scroll.\\
Let $l$ be a line in a fiber $F_\theta$ of $\ft$; since this contraction has
length $r$ we have
$$\det (\E_{|F_\theta})=(\det \E)_{|F_\theta}=(-K_Y)_{|F_\theta}=\Ol_{F_\theta}(r),$$ 
so the splitting type of $\E$ on $l$ is constantly $\Ol_\pu(1)^{\oplus r}$;
it follows that $\E_{|F_\theta} \simeq \Ol_{F_\theta}(1)^{\oplus r}$.
Therefore $p^{-1}(F_\theta)=\proj_{F_\theta}(\E_{|F_\theta}) 
\simeq \proj^{r-1}\times F_\theta$; since $p^{-1}(F_\theta)=\pst^{-1}(\pst(p^{-1}(F_\theta)))$
the subvarieties $\{x\} \times F_\theta$ of $\proj^{r-1} \times F_\theta$ correspond
to fibers of $\pst$.\\
In particular, if  $J \simeq \proj^r$ is a jumping fiber of $\ft$, then
$p^{-1}(J)= \proj_J(\E_{|J}) \simeq \proj^{r-1} \times J \simeq \proj^{r-1} \times \proj^r$ 
and the restriction $\pst:p^{-1}(J) \to \pst(p^{-1}(J))$ is a fibration in $\proj^r$,
hence each fiber is a jumping fiber.\qed 


\section{Recognizing products}\label{recprod}

In this section we collect some technical results that we are going to use
in order to establish whether a ruled Fano manifold is a product of another Fano manifold
with a suitable projective space.\par
\medskip
The idea of the following lemma is taken from \cite[Lemma 1.2.2]{AWI}.

\begin{lemma}\label{products}
Let $(X,Y,\E)$ be a $\proj^{r-1}$-ruled Fano manifold, and let $R_\E \subset \cone(X)$ 
be the extremal ray corresponding to the bundle projection.
Suppose that there exist an open subset $X^0 \subset X$ and a proper morphism 
$\psi:X^0 \to Z$ onto a variety $Z$ of dimension $r-1$ which does not contract curves of $R_\E$. 
Then $X \simeq \proj^{r-1} \times Y$.
\end{lemma}

{\bf Proof.} \quad Let $F$ be a general fiber of $\psi$; the dimension of $F$ is 
$\dim F= \dim X- \dim Z= \dim Y$, therefore $F$ dominates $Y$, since $\psi$ does not contract 
curves in the fibers of $p$.\\
Denote by $p_F:F \to Y$ the restriction of $p$ to $F$ and consider the 
pullback $\E_F=p_F^*\E$; denoted by $X_F$
the projectivization $\proj_F(\E_F)$, we have a commutative diagram
$$
\xymatrix@=40pt{X_F  \ar[r]^{\widetilde p_F} \ar[d]_{\widetilde p}  &X \ar[d]^{p}\\
F \ar[r]_{p_F}  &Y}
$$
By the universal property of the fiber product, $\widetilde p$ has a section
$s:F \to X_F$ such that $\widetilde p_F \circ s$ is the embedding of $F$ into $X$.
Let $\widetilde F=s(F)$ be the image of $F$ in $X_F$; by the canonical bundle
formula for $X_F$ we have
$$r\xi_{\E_F}-{\widetilde p}{\,}^*\det \E_F= 
-K_{X_F} + {\widetilde p}{\,}^*K_F.$$
Since ${\widetilde p}{\,}^*K_F
= K_{\widetilde F}= (K_{X_F})_{|\widetilde F}$, restricting to $\widetilde F$ we have
$(r\xi_{\E_F}-{\widetilde p}{\,}^*\det \E_F)_{|\widetilde F} =\Ol_{\widetilde F}$; 
therefore, using the canonical bundle formula for $X$,
$$\Ol_{F}=(r\xi_{\E}-{p}^*\det \E)_{|F} =(-K_{X} + {p}^*K_Y)_{|F}.$$
It follows that $\Ol_{F}= (K_X)_{|F} = p_F^*K_Y$, so $p_F$ is unramified. 
As $Y$, being Fano, is simply connected $p_F$ is an isomorphism,
hence $F$ is a section of $p$. To this section it
is associated an exact sequence of bundles over $Y$
\begin{equation}\label{seqY}
 \shse{\E'}{\E}{H}
\end{equation}
such that $F \in H^0(\xi_\E \otimes p^*\E'^\vee)$;
in particular the normal bundle of $F$ in $X$ is $(\xi_\E \otimes p^*\E'^\vee)_{|F}$.\\
Pulling back the sequence \ref{seqY} to $F$ we obtain an exact sequence of bundles
over $F$
\begin{equation}\label{seqF}
\shse{p_{F}^*\E'}{p_{F}^*\E}{p_{F}^*H}
\end{equation}
Since $F$ is a general fiber of $\psi$, its normal bundle in $X$ is trivial; thus we have 
$$\Ol_F^{\oplus r-1} = N_{F/X}= (\xi_\E \otimes p^*\E'^\vee)_{|F}.$$ 
It follows that $(p^*\E')_{|F} \simeq (\xi_\E)_{|F}^{\oplus r-1}$; 
therefore we can rewrite the sequence \ref{seqF} as
\begin{equation}\label{seqF2}
\shse{\xi_{\E_F}^{\oplus r-1}}{\E_F}{p_{F}^*H}.
\end{equation}
Recalling that
$(\det \E)_{|F}=r {\xi_\E}_{|F}=r\xi_{\E_F}$, we have $p_F^*H=\xi_{\E_F}$
and the sequence \ref{seqF2} splits, because $h^1(F,\Ol_F)=0$.
Thus $\E_F$ is decomposable as $\xi_{\E_F}^{\oplus r}$ and,
being $p_F$ is an isomorphism, also $\E$ is decomposable, as
a sum of $r$ copies of $H$. \qed \par

\begin{remark} In the proof of the lemma, instead of assuming that $Y$ is a Fano manifold,
it is enough to assume that $Y$ is simply connected and that $h^1(Y,\Ol_Y)=0$.
\end{remark}

\begin{corollary}\label{ruledoverprod}
Let $(X,Y,\E)$ be a $\proj^{r-1}$-ruled
Fano manifold of index $r$; assume that $Y \simeq \proj^{r-1} \times W$ and
denote by $\pi_1$ and $\pi_2$ the projections of $Y$ onto the factors. Then there exists
a vector bundle $\E'$ over $W$ such that $\pi_2^* \E'= \E \otimes \pi_1^*\Ol_{\proj^{r-1}}(-1)$ 
and $X = \proj^{r-1} \times \proj_W(\E')$.
\end{corollary}

{\bf Proof.} \quad The projection $\pi_2$ is the contraction associated to an extremal ray 
$\theta \subset \cone(Y)$; let $\vartheta \subset \cone(X)$ be its fellow ray.
By proposition \ref{spcontr} the contraction associated to $\vartheta$, 
$\psi_\vartheta: X \to Z$, is a $\proj^{r-1}$-bundle and $Z=\proj_W(\E')$, 
with $\E \otimes \pi_1^*\Ol_{\proj^{r-1}}(-1) = \pi_2^* \E'$.\\
In particular there exists a vector bundle $\F$ over $Z$ such that $(X,Z,\F)$
is a $\proj^{r-1}$-ruled Fano manifold; we can apply lemma \ref{products} to 
$(X,Z,\F)$, taking as $\psi$ the composition $\pi_1 \circ p:X \to \proj^{r-1}$.\qed\par

\begin{proposition}\label{rcprod} Let $(X,Y,\E)$ be a $\proj^{r-1}$-ruled
Fano manifold of index $r$. 
Suppose that there exist $R_1, \dots R_{\rho_Y}$ extremal  rays of length $r$ 
in $\cone(Y)$ such that $Y$ is rationally connected with respect
to curves in the corresponding families $R^1, \dots, R^{\rho_Y}$ 
(see example \ref{unex}). Then $X \simeq \proj^{r-1} \times Y$.
\end{proposition}

{\bf Proof.} \quad Let $C_i$ be a curve in the family $R^i$; since $\E$ is ample 
and $\det \E \cdot C_i=-K_Y \cdot C_i=l(R_i)=r$, denoting by $f_i:\pu \to C_i$ 
the normalization morphism, we have $f_i^*\E=\Ol_\pu(1)^{\oplus r}$.\\
Let $X_i= \pu \times_Y X= \proj_\pu(f_i^*\E) =\pu \times \proj^{r-1}$ and
let $G_i$ be the image of $X_i$ in $X$.\\ 
We have a commutative diagram
$$
\xymatrix@=40pt{X_i  \ar[r]^{\bar f_i} \ar[d]_{\bar p}  &X \ar[d]^{p}\\
\pu \ar[r]_{f_i}  &Y}
$$
Let ${\widetilde C}_i$ be a section of $\bar p:X_i \to \pu$, let 
$\Gamma_i=\bar f_i({\widetilde C}_i)$ be its image
 in $X$ and let $V^i$ be a family of deformations of $\Gamma_i$; 
by the canonical bundle formula we have $-K_X \cdot \Gamma_i = r\xi_\E \cdot \Gamma_i =r$, 
therefore the family $V^i$ is an unsplit family.\\
Let $x$ be a point of $X$ and $y$ a point of $Y$; as $Y$ is rationally connected 
with respect to curves in $R^1, \dots, R^{\rho_Y}$, there exists a chain of curves 
$C_{i_1}, \dots, C_{i_m}$ in $R^1, \dots, R^{\rho_Y}$ connecting $p(x)$ and $y$, with
$m \le 2^{\dim Y} -1$.\\
Let $y_1$ be a point in $C_{i_1} \cap C_{i_2}$ and let
$\Gamma_{i_1}$ be a curve in $V^{i_1}$ which is mapped to $C_{i_1}$ and passes
through $x$. 
The fiber of $p$ over $y_1$ is contained in $G_{i_2}$, so
there is a minimal section $\Gamma_{i_2}$ in $G_{i_2}$ which meets $\Gamma_{i_1}$;
repeating the argument we construct a chain of curves in $V^1, \dots, V^{\rho_Y}$
which joins $x$ with a point of the fiber over $y$.
We have thereby proved that, for every $x \in X$ and for some $m$, 
$\cloc_m(V^1, \dots, V^{\rho_Y})_{x}$ dominates $Y$.\par
\smallskip
Let $\psi:X^0 \to Z$ be the rc$(V^1, \dots, V^{\rho_Y})$-fibration; a general
fiber $F$ of $\psi$ is an equivalence class for the rc$(V^1, \dots, V^{\rho_Y})$-relation, 
thus it contains $\cloc_m(V^1, \dots, V^{\rho_Y})_{x}$ for every point $x \in F$ 
and every $m$; then we have $\dim F \ge \dim Y$ and $\dim Z \le \dim X -\dim F \le r-1$. \\
On the other hand, $F$ cannot contain a curve in a fiber of $p$,
otherwise $R_\E$ would be contained in the subvector space of $N_1(X)$ 
generated by the classes of ${V}^1,\dots, {V}^{\rho_Y}$ by \cite[Corollary 4.2]{ACO}.
Being $\loc(R^\E)_F=X$, this, again by \cite[Corollary 4.2]{ACO}, would imply that 
the class of every curve in $X$ would be contained in the subvector space of $N_1(X)$ 
generated by the classes of ${V}^1,\dots, {V}^{\rho_Y}$, hence
$\rho_X=\rho_Y$, a contradiction.\\
In particular it follows that $\dim F= \dim Y$; therefore $\dim Z=r-1$ and we can 
apply lemma \ref{products} to $(X,Y,\E)$ and $\psi$ to conclude.\qed \par


\section{Fano manifolds with many fiber type contractions}\label{mafi}

In this section we will prove that a ruled Fano fivefold of index
two and Picard number greater than three is a product.
We will derive this conclusion from two more general results
concerning Fano manifolds with many fibrations. 

\begin{proposition}\label{fibertype} Let $X$ be a Fano manifold 
of dimension $n$ and pseudoindex $i_X \ge 2$ which has only contractions of fiber type. 
Then $\rho_X \le n$. Moreover,
\begin{enumerate}
\item if $\rho_X=n$, then $X=(\pu)^n$;
\item if $\rho_X=n-1$, then $X=(\pu)^{n-2} \times \pd$ or $X=(\pu)^{n-3} 
\times \proj_{\pd}(T\pd)$.
\end{enumerate}
\end{proposition}

{\bf Proof.} \quad By \cite[Theorem 2.2]{Wi1} we have that a Fano manifold
of dimension $n$ admits at most $n$ fiber type elementary contractions, and the bound on the
Picard number follows. More precisely we have that the cone of curves of $X$
is generated by at most $n$ extremal rays.\par
\smallskip 
We can assume that $n \ge 4$, since for lower dimensions the claimed
result follows from the classification of Fano manifolds.\\
Suppose that $\rho_X=n$; by the discussion above we have
$\cone(X)=\langle R_1, \dots, R_n \rangle$.
Let $R^1, \dots, R^n$ be the corresponding families of rational curves, as in example
\ref{unex}; by 
\cite[Lemma 5.4 (c)]{ACO} we have 
$$n \ge \dim \loc(R^1, \dots, R^n)_x \ge \sum_{i=1}^n (-K_X \cdot R^i-1) \ge n,$$
forcing $-K_X \cdot R^i=2$ for every $i$ (recall that $i_X \ge 2$) 
and $\sum_{i=1}^n (-K_X \cdot R^i-1)=n$. We can therefore apply
\cite[Theorem 1]{Op} to conclude.\par
\smallskip
Suppose now that $\rho_X=n-1$; let $R_1, \dots, R_{n-1}$
be extremal rays of $X$ which span $N_1(X)$ and let $R^1, \dots, R^{n-1}$
be the corresponding families of rational curves.\\
Suppose that, among the chosen rays, there exists a ray $R_{i(1)}$ 
such that the associated contraction $\f_{i(1)}$ has a fiber 
$F$ of dimension greater than one.
We claim that for every ray $R_{i(j)} \in \{R_1, \dots, R_{n-1}\}$ different
from $R_{i(1)}$ the contraction associated to $R_{i(j)}$ 
is equidimensional with one dimensional fibers.\\
Assume by contradiction that there exists an index $i(2)$ 
such that the contraction associated to $R_{i(2)}$ has a fiber $G$ of
dimension $\ge 2$.\\
Consider an irreducible component $D$ of $\loc(R^{i(3)}, \dots, R^{i(n-1)})_{G}$, 
which, by \cite[Lemma 5.4 (c)]{ACO}, has dimension 
$$ \dim D \ge \sum_{j=3}^{n-1} (-K_X \cdot R^{i(j)}-1) + \dim G \ge n-1.$$
By \cite[Lemma 5.1]{ACO}, $N_1(D)=\langle R_{i(2)}, \dots, R_{i(n-1)}\rangle$, therefore
we cannot have $D=X$, thus $D$ is an effective divisor in $X$. 
We will now derive a contradiction by considering the intersection number of this divisor
with the family $R^{i(1)}$.\\
Suppose first that $D \cdot R^{i(1)}>0$; in this case $D$ meets $F$,
which has dimension at least two, whence the intersection $D \cap F$ contains
a curve, contradicting the fact that curves in $R^{i(1)}$ are numerically independent
from curves in $D$.\\
Suppose now that $D \cdot R^{i(1)} = 0$ and let $C_{i(1)}$ be a curve of $R^{i(1)}$
meeting $D$. Since the intersection number is zero, this curve is contained in
$D$, contradicting again the independence of curves in $R^{i(1)}$ from curves in $D$.\par
\smallskip
We have thereby proved that $X$ has at least $n-2$ extremal rays whose associated
contractions are equidimensional with one dimensional fibers.
Let $\f_j:X \to Y_j$ be one of these contractions; by lemma \ref{pubundle} 
there exists a rank two vector bundle $\E_j$ on $Y_j$ such that $X= \proj_{Y_j}(\E_j)$.\\ 
By lemma \ref{bundlebase} (a), $Y_j$ is a Fano manifold of pseudoindex $i_{Y_j} \ge i_X \ge 2$ 
and, by part (c) of the same lemma, has only contractions of fiber type,
so, by induction on the dimension, $Y_j \simeq (\pu)^{n-3} \times \pd$ 
or $Y_j \simeq (\pu)^{n-4} \times \proj_{\pd}(T\pd)$.\\
It follows that $i_{Y_j}=2 =i_X$, hence, by lemma \ref{bundlebase} (b), the restriction of 
$\E_j$ to every fiber of a $\pu$-bundle contraction of $Y_j$ splits 
as a sum of two line bundles of the same degree.\\
Up to twist $\E_j$ with a suitable line bundle in $\pic(Y_j)$, we can now
assume that the restriction of $\E_j$ to any fiber of a $\pu$-bundle contraction is 
$\Ol_\pu(1) \oplus \Ol_\pu(1)$.\\ 
In particular $K_{Y_j}+ \det \E_j$ is trivial on all the
extremal rays of $Y_j$, hence $K_{Y_j} + \det \E_j=\Ol_{Y_j}$; by the canonical bundle
formula we have $-K_X=2\xi_{\E_j}$, consequently $(X,Y_j,\E_j)$ is a ruled
Fano manifold of index two.\\
For both possible basis $Y_j$ the ruled Fano manifold $(X,Y_j,\E_j)$
verifies the assumptions of proposition \ref{rcprod}, so we have $X= \pu \times Y_j$.\qed

\begin{proposition}\label{mixedtype} Let $X$ be a Fano manifold of dimension $n$
and pseudoindex $i_X \ge 2$ such that all its elementary contractions but one 
are of fiber type. Then $\rho_X \le n-1$, equality holding if and only if
$X=(\pu)^{n-3} \times Bl_p(\pt)$.
\end{proposition}

{\bf Proof.} \quad We can assume that $n \ge 4$, since for lower dimensions the claimed
result follows from the classification of Fano manifolds.\\
Let $R_1$ be the birational ray and let $R_2, \dots, R_{\rho_X}$ be
fiber type rays such that $R_1, R_2,\dots,$ $R_{\rho_X}$ span $N_1(X)$.
Let $\f_{1}:X \to X'$ be the contraction of $R_1$ and let $F$ be a nontrivial 
fiber of $\f_{1}$;  since $\f_1$ is birational, by proposition \ref{fiberlocus}
we have $\dim F \ge 2$.\\
For every permutation $i(2), \dots, i(\rho_X)$ of the
integers $2, \dots, \rho_X$, by \cite[Lemma 5.4 (c)]{ACO} we have
$$\dim \loc(R^{i(2)}, \dots, R^{i(\rho_X)})_{F} \ge \dim F+ \rho_X-1,$$
forcing $\rho_X \le n-1$; moreover, if equality holds, we have $\dim F=2$ and 
$X= \loc(R^{i(2)}, \dots, R^{i(\rho_X)})_{F}$.\\
In particular we note for later use that, since $\f_1$ is birational and all 
its nontrivial fibers have dimension $= 2$, $\Exc(\f_1)$ is a divisor 
by proposition \ref{fiberlocus}.\\
Set $T_{i(2)}=\loc(R^{i(2)})_F$; being 
$X= \loc(R^{i(3)}, \dots, R^{i(\rho_X)})_{T_{i(2)}}$, by \cite[Lemma 1]{Op}
every curve $C \subset X$ is equivalent to a linear combination 
$$\alpha \Gamma_{i(2)}+ \sum_{k = 3}^{\rho_X} \alpha_k R^{i(k)}$$
of a curve $\Gamma_{i(2)}$
in $T_{i(2)}$ and curves in $R^{i(3)}, \dots, R^{i(\rho_X)}$ with
$\alpha \ge 0$.
By \cite[Corollary 2.23]{CO} every curve in $T_{i(2)}$
is numerically equivalent (in $X$) to a linear combination with positive
coefficients of a curve in $F$ (and so whose numerical class is in $R_1$) 
and a curve in $R^{i(2)}$; hence we can write $C$ as a combination
$$\alpha_1 R_1 + \alpha_2 R^{i(2)} + \sum_{k =3}^{\rho_X} \alpha_k R^{i(k)}, $$
with $\alpha_1, \alpha_2 \ge 0$.\\
Since this is true for every permutation $i(2), \dots, i(\rho_X)$,
and the decomposition of $[C]$ is unique, we get that
$\alpha_{k} \ge 0$ for all $k$ and 
$\cone(X)=\langle R_1, R_2, \ldots, R_{\rho_X} \rangle$.\par
\medskip
Denote again by $F$ a nontrivial fiber of $\f_1$ and, for every $i=2, \dots, \rho_X$ 
consider an irreducible component $D_i$ 
of $\loc(R^{2}, \dots, \hat{R^i}, \dots, R^{\rho_X})_{F}$
which, by \cite[Lemma 5.4 (c)]{ACO}, has dimension 
$$ \dim D_i \ge \sum_{j=2}^{\rho_X} (-K_X \cdot R^j-1) +(K_X \cdot R_i +1) + \dim F \ge n-1. $$
By \cite[Lemma 5.1]{ACO}  $N_1(D_i)=\langle R_1, \dots \hat{R_i}, \dots, R_{\rho_X}\rangle$, 
therefore we cannot have $D_i=X$, whence $D_i$ is an effective divisor in $X$. \\
As in proposition \ref{fibertype} we can now prove that the contraction
$\f_{i}:X \to Y_i$, associated to the ray $R_{i}$, has one dimensional fibers,
since the intersection of this fibers with $D_i$ must be $0$-dimensional, hence, by 
lemma \ref{pubundle} there exists a rank two vector bundle $\E_i$ on $Y_i$ such that 
$X= \proj_{Y_i}(\E_i)$.\\ 
By lemma \ref{diveff}, for at least one index $j \in \{2, \dots,\rho_X\}$
we have $\Exc(\f_1) \cdot R_j >0$; let $\f_{j}:X \to Y_j$ be the contraction associated to 
the ray $R_j$.\\
By lemma \ref{bundlebase}, $Y_j$ is a Fano manifold of pseudoindex $i_{Y_j} \ge 2$;
by lemma \ref{fiberdim} all the extremal contractions of $Y_j$ are
of fiber type and, by the same lemma, one of these contractions has two dimensional fibers.
We can apply proposition \ref{fibertype} to $Y_j$ to get 
$Y_j \simeq (\pu)^{n-3} \times \pd$.\\
Let $p_1:Y_j \to \pu$ be the projection onto the first factor;
the projection to the other factors is an extremal elementary contraction 
$\ft:Y_j \to (\pu)^{n-4} \times \pd$, associated to a ray $\theta \subset \cone(Y_j)$.\\
Let $\vartheta \subset \cone(X)$ be the fellow ray of $\theta$;
since $\ft$ has one dimensional fibers, the same is true for
the contraction associated to $\vartheta$, $\pst:X \to Z$.
Therefore $\vartheta \not = R_1$, and the associated contraction $\pst$
is a $\pu$-bundle over a smooth Fano variety $Z$, which has pseudoindex
$i_Z \ge i_X \ge 2$.\\
Consider the following diagram
$$
\xymatrix@C=30pt@R=40pt{X  \ar[d]_{\pst} \ar[r]^{\f_j}  & Y_j \ar[d]^{\ft} \ar[r]^{p_1} & \pu\\
Z \ar[r] \ar[r]  &(\pu)^{n-4} \times \pd }
$$

We can apply lemma \ref{products} to $X$ and $\psi=p_1 \circ \f_j: X \to \pu$ 
and obtain $X \simeq \pu \times Z$.
It follows that $Z$ has a birational contraction, so, by induction
$Z \simeq (\pu)^{n-4} \times Bl_p(\pt)$ and $X \simeq (\pu)^{n-3} \times Bl_p(\pt)$.
\qed \par

\begin{corollary} Let $X$ be a Fano fivefold of index $r_X \ge 2$ and Picard number
$\rho_X \ge 4$. Then
\begin{enumerate}
\item $X \simeq (\pu)^5$;
\item $X \simeq (\pu)^2 \times \proj_\pd(T\pd)$;
\item $X \simeq (\pu)^2 \times Bl_p(\pt)$.
\end{enumerate}
\end{corollary}

{\bf Proof.} \quad Note that, since $\rho_X \ge 4$, we have $i_X \le 2$, by
\cite[Theorem 1.4]{ACO}, hence $r_X=i_X=2$.\\
By \cite[Theorem 1.1]{CO}, if $\rho_X \ge 4$, then $X$ has at most
one birational contraction, and the conclusion follows from propositions
\ref{fibertype} and \ref{mixedtype}. Note that $(\pu)^{n-2} \times \pd$
has been excluded since its index is one. \qed \par


\section{Proof of theorem \ref{rho3}: Classification of the base}

In this section we begin the study of ruled Fano fivefolds $(X,Y,\E)$ of index
two and Picard number three, which is the most complicated case.\\
We start by considering the possible bases $Y$ such that
there exists a ruled Fano fivefold $(X,Y,\E)$ as above which is not a product.
By lemma \ref{bundlebase} $Y$ is a Fano fourfold of pseudoindex $i_Y \ge 2$,
and $\rho_Y=2$, since we are assuming $\rho_X=3$.
We will give a complete classification of fourfolds $Y$ as above
which have a birational contraction (Proposition \ref{Ybir}), 
and a more rough one of the ones with two fiber type contractions
(Proposition \ref{Yfib}). Then, using the criteria for recognizing products
previously estabilished, we will show that there are only four possibilities
for $Y$  (Proposition \ref{noprod}). 

\begin{proposition}\label{Ybir} Let $Y$ be a Fano fourfold of pseudoindex $i_Y \ge 2$ 
and Picard number $\rho_Y=2$ such that the contraction $\f_\theta:Y \to Y'$,
associated to one extremal ray $\theta \subset \cone(Y)$, is birational.  
Then $Y$ is one of the following:
\begin{enumerate}
\item $Bl_p(\proj^4)$ with $p$ a point in $\proj^4$;
\item $Bl_l(\proj^4)$ with $l$ a line in $\proj^4$;
\item $Bl_l(\quadr^4)$ with $l$ a line in $\quadr^4$;
\item $Bl_\Gamma(\quadr^4)$ with $\Gamma$ a conic in $\quadr^4$ not 
contained in a plane $\Pi \subset \quadr^4$;
\item $\proj_\pt(\Ol_\pt \oplus \Ol_\pt(2))$; 
\item $\proj_\qt(\Ol_\qt \oplus \Ol_\qt(1))$.
\end{enumerate}
\end{proposition}

{\bf Proof.} \quad The cone of curves of $Y$ is generated by two extremal rays:
$\cone(Y)= \langle \theta, \oth \rangle$.\\
The length of every extremal ray on a Fano manifold is clearly greater than or
equal to the pseudoindex; moreover, for a birational extremal ray, by proposition
\ref{fiberlocus}, the length is bounded above by the dimension of the manifold minus one,
hence
$$2 \le l(\theta) \le 3.$$

If $l(\theta)=3$, by \cite[Theorem 1.1]{AO}, the associated contraction
$\f_\theta:Y \to Y'$ is the blow up at a point of a smooth
variety $Y'$; Fano manifolds which are the blow up at a point of a smooth variety 
are classified in \cite[Theorem 1.1]{BCW}, which gives three possible cases. 
Among these cases only the blow up at a point of the projective space has pseudoindex greater
than one, hence we are in case (1).\par
\smallskip
If $l(\theta)=2$, by \cite[Theorem 5.2]{AO}, either $\f_\theta$ is the blow up 
of a smooth variety along a smooth curve, or its exceptional locus $\Exc(\f_\theta)$ 
is isomorphic to $\pt$ or to a (possibly singular) three dimensional quadric 
and $\f_\theta(\Exc(\f_\theta))$ 
is a point.\par
\smallskip
If $\f_\theta:Y \to Y'$ is the blow up of a smooth variety along a smooth curve, we can apply 
\cite[Theorem 1.3]{AO4} and, recalling that we are assuming $\rho_Y=2$, 
we have cases (2), (3) and (4).\par 
\smallskip
If else $\f_\theta(\Exc(\f_\theta))$ is a point, we consider the contraction 
$\f_{\oth}:Y \to T$, associated to the extremal ray $\oth$; the effective divisor 
$\Exc(\f_\theta)$ is positive on $\oth$ by lemma \ref{diveff}, therefore, by corollary 
\ref{onedim}, $\f_{\oth}$ makes $Y$ a $\pu$-bundle over  $T$, $Y=\proj_T(\F)$.
We can thus apply lemma \ref{fano3fold}, obtaining that either $Y$ is a product, or $T$ 
is a projective space or a smooth quadric. The first case has to be excluded since 
$\pu \times T$ does not have a birational contraction; in the second case we note that 
$\F$ is a Fano bundle on $T$, whence we can use the classification 
in \cite{SW2}, looking for bundles such that their projectivization has pseudoindex $\ge 2$ 
and a birational extremal contraction.\\
By that classification it turns out that the only possibilities are 
number (5) and (6) in our list.\qed\par

\begin{proposition}\label{Yfib} Let $Y$ be a Fano fourfold of pseudoindex $i_Y \ge 2$ 
and Picard number $\rho_Y=2$ with two fiber type extremal contractions. Then $Y$ is one of the
following:
\begin{enumerate}
\item[(1)] a product $\pu \times W$;
\item[(2)] a variety whose extremal rays have length $2$ and associated contractions
with fibers of dimension $\le 2$;
\item[(3)] $\pd \times \pd$;
\item[(4)] $\proj_\pd(T\pd(-1) \oplus \Ol_\pd)$.
\end{enumerate}
\end{proposition}

{\bf Proof.} \quad The manifold $Y$ is Fano and has Picard number two, 
so its cone is spanned by two extremal rays: 
$\cone(Y)=\langle \theta,\overline{\theta} \rangle$.\par
\medskip
Suppose that the contraction associated to one extremal ray, say $\oth$,
has a three dimensional fiber $F_{\,\oth}$; then, by lemma \ref{diveff}, 
$F_{\,\oth} \cdot \theta >0$. By corollary \ref{onedim}, the contraction of
 $\theta$, $\f_\theta:Y \to W$, 
makes $Y$ into a $\pu$-bundle over a smooth threefold
$W$, $Y = \proj_W(\F)$;
by lemma \ref{fano3fold}, either $Y \simeq \pu \times W$, or
$W$ is $\proj^3$ or $\quadr^3$.\\
By the classification given in \cite{SW2}, 
there are no of Fano bundles over $\pt$ and $\quadr^3$ such that
their projectivization is not a product and has two fiber type contractions,
one of which has a three dimensional fiber.\par
\medskip
Therefore either we are in case (1) or both the contractions of $Y$ have fibers 
of dimension $\le 2$; this implies that the lengths of the extremal rays are $\le 3$, by
proposition \ref{fiberlocus}.\\
Either we are in case (2) or the length of one extremal ray, say $\theta$, is equal to
three; again by proposition \ref{fiberlocus} we have that
$\f_\theta:Y \to W$ is equidimensional with fibers
of dimension two.\\ 
By lemma \ref{brauer} $W$ is smooth and so, being a smooth surface of 
Picard number one dominated by a Fano manifold, $W \simeq \pd$; moreover,
by the same lemma $Y=\proj_\pd(\F)$ for some rank three vector 
bundle on $\pd$. In particular $\F$ is a Fano bundle over $\pd$.\\ 
From the classification of such bundles given in \cite{SW1}, recalling that, in our case,
the other contraction of $Y$ has length $\ge 2$, we are either in case (3)
or in case (4).\qed\par

\begin{proposition}\label{noprod} Suppose that there exists a ruled Fano fivefold of index two 
$(X, Y, \E)$ with $\rho_X=3$ which is not a product with $\pu$ as a factor.
Then $Y$ is one of the following:
\begin{enumerate}
\item $Bl_p(\proj^4)$;
\item $Bl_l(\proj^4)$;
\item $\pd \times \pd$;
\item $\proj_\pd(T\pd(-1) \oplus \Ol_\pd)$.
\end{enumerate}
\end{proposition}

{\bf Proof.} \quad Suppose first that $Y$ has a birational contraction; then
$Y$ is one of the manifolds listed in proposition \ref{Ybir}.
The varieties (3)-(6) are rationally connected with respect to minimal curves 
in the extremal rays, which have length two, so, if they are the base
of a ruled Fano fivefold $(X, Y, \E)$ of index two, then $X$ is a product 
$\pu \times Y$ by proposition \ref{rcprod}.
Therefore, if $Y$ has a birational contraction and $X$ is not a product,
$Y$ is either $Bl_p(\proj^4)$ or $Bl_l(\proj^4)$ (cases (1) and (2) of
proposition \ref{Ybir}).\par
\medskip 
Suppose now that $Y$ has only fiber type contractions; then, by proposition
\ref{Yfib}, we have four possible cases. To finish the proof we have
to rule out cases (1) and (2) of that proposition.\par
If $Y \simeq \pu \times W$, we can apply corollary \ref{ruledoverprod} to get that $X$
is a product $\pu \times \proj_W(\E')$.\par
\medskip
We are left with the case of a manifold $Y$ whose extremal rays have length $2$ 
and associated contractions with fibers of dimension $\le 2$.
Let $\theta$ be one of the rays in $\cone(Y)$, let $\f_\theta:Y \to W$
be the associated contraction and let $R^\theta$ be the associated family 
of rational curves; we claim that $R^\theta$ is a covering family.\\
If the general fiber $F_\theta$ of $\f_\theta$ has dimension one, 
this follows from proposition \ref{iowifam}, since $\loc(R^\theta)_x$ is
contained in the fiber of $\f_\theta$ through $x$:
$$\dim \loc(R^\theta) \ge \dim Y +l(\theta) - 1 - \dim \loc(R^\theta)_x \ge 4.$$
If else $F_\theta$ has dimension two, then, by adjunction,
it is a smooth quadric and therefore it is covered by curves in $R^\theta$,
which are lines in the quadric.\\
We can thus consider the rc$(R^\theta,R^{\oth})$-fibration,
whose image has to be a point, being $\rho_Y=2$.
It follows that $Y$ is rationally connected with respect
to curves in $R^\theta$ and $R^{\oth}$ and $X$ is a product $\pu \times Y$
by proposition \ref{rcprod}.\qed \par


\section{Proof of theorem \ref{rho3}}

In this section we achieve the classification of ruled fivefolds $(X,Y,\E)$ of index two 
and Picard number three, proving theorem \ref{rho3}.\par
\medskip
First we prove that, if $X$ is not a product, one of the contractions
of $X$ is birational (proposition \ref{nofib}). We then consider separately 
the case in which also $Y$ has a birational contraction
(proposition \ref{XYbir}) and the case in which both the contractions of $Y$ 
are of fiber type (proposition \ref{XbirYnot}).

\setcounter{equation}{0}
\begin{proposition}\label{nofib} Let $(X,Y,\E)$ be a ruled Fano fivefold of index two 
with $\rho_X = 3$ such that $X$ has only fiber type contractions. Then $X$ is a product
with $\pu$ as a factor.
\end{proposition}

{\bf Proof.} \quad Since $X$ has only fiber type contractions, the same is true
also for $Y$ by corollary \ref{excloc}, so, by proposition \ref{noprod}, 
if $X$ is not a product with $\pu$ as a factor, then $Y$ is either $\pd \times \pd$ or 
$\proj_\pd(T\pd(-1) \oplus \Ol_\pd)$.\par
\medskip
{\bf Case a)} \quad $Y \simeq \pd \times \pd$.\par
\medskip
The cone of curves of $Y$ is generated by two extremal rays, $\theta$ and 
$\bar \theta$, corresponding to the projections $\f_{\theta}, \fot:Y \to \pd$
Let $\vartheta$ and $\ovth$ be the fellow rays of $\theta$ and $\oth$, respectively,
and denote by $\pst:X \to Z$ and $\psot:X \to \overline{Z}$ the associated
contractions.
By proposition \ref{spcontr2} the contractions $\pst$ and $\psot$
are $\pu$-bundles and $p^*\Ol_Y(1,1)$ restricts to $\Ol_\pu(1)$
on the fibers of $\pst$ and $\psot$.
Hence there exist two vector bundles
$\F$ on $Z$ and $\overline{\F}$ on $\overline{Z}$
such that $(X,Z,\F)$ and $(X,\overline{Z},\overline{\F})$ are ruled Fano 
fivefolds of index two.\\
Since all the contractions of $X$ are of fiber type, the same is true
also for $Z$ and $\overline{Z}$, by corollary \ref{excloc}.
We can apply proposition \ref{noprod} to $(X,Z,\F)$ and to $(X,\overline{Z},\overline{\F})$
and we have for $Z$ and $\overline{Z}$ two possibilities: $\pd \times \pd$
or $\proj_\pd(T\pd(-1) \oplus \Ol_\pd)$.\\
In the last case one extremal
contraction of $Z$ ($\overline{Z}$) is a special \ba scroll onto $\pt$ so, by proposition
\ref{spcontr3}, also one contraction of $X$ has to be a special \ba scroll with
jumping fibers, but we have already proved that all the contractions of $X$ are $\pu$-bundles.\\
It follows that both $Z$ and $\overline{Z}$ are $\pd \times \pd$.
All the extremal rays of $X$ have length two, hence $\xi_\E$ restricts
to $\Ol_\pu(1)$ on the fibers of any contraction of $X$.\\
Consider the commutative diagram
$$
\xymatrix@=15pt{\pd& & & & Z \ar[rd] & \\
& Y \ar[lu]_{\ft} \ar[ld]^{\fot}& & 
X \ar[ll]^-{p} \ar[ru]^{\pst} \ar[rd]_-{\psot} \ar[rr]^-{\psi_\sigma} & & \pd \\
 \pd & & & & \overline{Z} \ar[ru] & }
$$

The line bundle $\xi_{\E} \otimes p^*\Ol_Y(-1,-1)$ is trivial on the face $\sigma$
spanned by $\vartheta$ and $\ovth$, and restricts to $\Ol_\pu(1)$ on the fibers of
$p$, hence $\xi_{\E(-1,-1)}=\xi_{\E} \otimes p^*\Ol_Y(-1,-1)=\psi_{\sigma}^* \Ol_\pd(1)$
is spanned. Equivalently $\E(-1,-1)$ is spanned and $h^0(\E(-1,-1))=3$.
We thus have a surjective map $\Ol_{Y}^{\oplus 3} \to \E(-1,-1) \to 0$, which gives rise
to an exact sequence
$$\shse{L}{\Ol_{Y}^{\oplus 3}}{\E(-1,-1)};$$
computing the splitting type we find $L \simeq \Ol_{Y}(-1,-1)$.
The dual bundle $L^\vee$ is thereby ample, therefore, by \cite[12.1.6]{Ful}, 
the map $L \to \Ol_Y^{\oplus 3}$ must have a non empty degeneracy locus,
whence $X=\proj_Y(\E(-1,-1)) \hookrightarrow \proj_Y(\Ol_Y^{\oplus 3})$
is not a $\pu$-bundle over $Y$, a contradiction.\par
\medskip
{\bf Case b)} \quad $Y \simeq \proj_\pd(T\pd(-1) \oplus \Ol_\pd)$.\par
\medskip
The cone of curves of $Y$ is generated by two extremal rays: $\theta$,
corresponding to the projection $\f_{\theta}:Y \to \pd$,
and $\oth$, corresponding to the contraction $\fot:Y \to \pt$, which
is a special B\v anic\v a scroll with exactly one jumping fiber $J\simeq \pd$,
which is the section corresponding to the trivial summand of the bundle
$T\pd(-1) \oplus \Ol_\pd$.\\
Let $\vartheta$ and $\ovth$ be the fellow rays of $\theta$ and $\oth$, respectively,
and denote by $\pst:X \to Z$ and $\psot:X \to \overline{Z}$ the associated
contractions.
By proposition \ref{spcontr2} the contraction $\pst:X \to Z$ is
a $\pu$-bundle, while, by proposition \ref{spcontr3}, the contraction
$\psot:X \to \overline{Z}$ is a special B\v anic\v a scroll with
a one parameter family of jumping fibers which are sections of $p$
over over $J$.\\
Since $\pst:X \to Z$ is a $\pu$-bundle, there exists a vector bundle
$\F$ on $Z$ such that $(X,Z,\F)$ is a ruled Fano fivefold of index two.
All the contractions of $Z$ are of fiber type by corollary \ref{excloc},
so proposition \ref{noprod} applied to $(X,Z,\F)$ gives us two possibilities: 
either $Z \simeq \pd \times \pd$ or $Z \simeq \proj_\pd(T\pd(-1) \oplus \Ol_\pd)$.
In the first case we conclude as in case a), replacing $(X,Y,\E)$ with
$(X,Z,\F)$, otherwise we consider the following commutative diagram
\objectmargin+{4pt}
$$
\xymatrix@=30pt{
& J \ar@{^{(}->}[d]  & p^{-1}(J) \ar@{^{(}->}[d] &\\
\pd & Y \ar[l]^{\ft} \ar[d]_-{\fot} &  X \ar[l]_-{p} \ar[r]^{\pst} \ar[d]_-{\psot}  & Z \ar[d]\\
&  \pt &  \overline{Z} \ar[l]^-{p'} \ar[r]  & \pt}
$$
By proposition \ref{spcontr3} there is an
isomorphism $f:\pu \times \pd \to p^{-1}(J)$, and the subsets $f(\{x\} \times \pd)$
are jumping fibers of $\psot$.
In particular the numerical class of every curve in $p^{-1}(J)$ 
belongs to the face $\langle{R_\E,\ovth}\rangle$.
It follows that $\pst$ is finite to one on $p^{-1}(J)$, but this is a contradiction 
since, by lemma \ref{fiberdim} every jumping fiber of $\psot$ has to be mapped by $\pst$ 
to a jumping fiber of the contraction $Z \to \pt$, but this map has only one 
jumping fiber.\qed\par

\begin{proposition}\label{XYbir} Let $(X,Y, \E)$ be a ruled Fano fivefold of index two 
with $\rho_X = 3$ such that both $X$ and $Y$ have a birational contraction.
Then, if $X$ is not a product with $\pu$ as a factor, one of the following happens:
\begin{enumerate}
\item $X \simeq Bl_p(\proj^4) \times_\pt Bl_p(\proj^4)$;
\item $X \simeq Bl_S(Bl_p(\proj^5))$ with $S$ the strict trasform of a plane $\ni p$.
\end{enumerate}
In these cases the corresponding pairs $(Y,\E)$ are, respectively,
\begin{enumerate}
\item $(Bl_p(\proj^4), 2H +E \oplus 3H +E)$, 
 $E$  exceptional divisor and $H$ pullback on $Y$ of $\Ol_\pt(1)$;
\item $(Bl_l(\proj^4), 2H-E \oplus 3H-E)$, 
$E$  exceptional divisor and $H$ pullback on $Y$ of $\Ol_{\proj^4}(1)$.
\end{enumerate}
\end{proposition}

{\bf Proof.} \quad We assume that $X$ is not a product and that $Y$ has a birational
contraction, so, by corollary \ref{noprod}, $Y$ is the blow up of $\proj^4$ 
 either along a point or along a line. \par
\medskip
{\bf Case a)} \quad $Y =Bl_p(\proj^4)$.\par
\medskip
Another possible description of $Y$ is $\proj_\pt(\Ol_\pt \oplus \Ol_\pt(-1))$;
let $\theta \subset \cone(Y)$ be the extremal ray corresponding to the 
$\pu$-bundle contraction $\ft:Y \to \proj^3$, let $E$ be the exceptional $\proj^3$
and let $H$ be the pullback of $\Ol_\pt(1)$.
Let $\vartheta \subset \cone(X)$ be the fellow ray of $\theta$;
by proposition \ref{spcontr}, the contraction 
associated to $\vartheta$, $\pst:X \to Z$, is a $\pu$-bundle, too. 
Moreover, by the same proposition, since $E$ restricts to $\Ol_{\proj^1}(1)$
on the fibers of $\ft$, we have $\E \otimes (-E)=\ft^*\E'$ and $Z=\proj_\pt(\E')$.\par
\smallskip 
Since $E_{|E} \simeq \Ol_\pt(-1)$ and $E$ is a section of $\ft$, we have
$$\E_{|E}= (\f_\theta^*\E' \otimes E)_{|E} \simeq \E'(-1).$$
Recalling that $(\det \E)_{|E}=(-K_Y)_{|E}=\Ol_\pt(3)$ and that $\E$ is ample,
we see that the splitting type of $\E$ on lines of $E$ is constantly
$\Ol_\pu(1) \oplus \Ol_\pu(2)$, hence, by \cite[Theorem 3.2.3]{okonek}, $\E_{|E}$ 
is decomposable as $\E_{|E} \simeq \Ol_\pt(1) \oplus \Ol_\pt(2)$.
It follows that $\E'\simeq \Ol_\pt(2) \oplus \Ol_\pt(3)$, thus
$\E \simeq (2H \oplus 3H) \otimes E$.\par
\medskip
{\bf Case b)} \quad $Y =Bl_l(\proj^4)$.\par
\medskip
Let $\theta \subset \cone(Y)$ be the extremal ray whose associated contraction, 
$\ft:Y \to \proj^4$, is the blow up of $\proj^4$ along a line.
Denote by $E$ the exceptional locus of $\ft$
and by $H$ the pullback of the ample generator of $\pic(\proj^4)$.\\
Let $\vartheta \subset \cone(X)$ be the fellow ray of $\theta$;
by proposition \ref{spcontr}, the associated contraction, 
$\pst:X \to X'$, is the blow up of a smooth fivefold along a smooth surface.\\
By the same proposition, since $-E$ restricts to $\Ol_\pd(1)$
on the fibers of $\ft$, there exists a rank two vector bundle on $X'$ such that
$\E \otimes E=\ft^*\E'$ and $X'=\proj_W(\E')$; by \cite[Lemma 2.10]{AO2} $\E'$ is ample.\\
The canonical bundle formula for blow ups, $K_{Y} = \ft^*K_{\proj^4}+2 E$,
combined with the determinant formula, $\det \ft^* \E' = \det \E +2 E$, gives
$$\ft^*(K_{\proj^4}+\det \E')= K_Y +\det \E=\Ol_Y,$$
whence $K_{\proj^4}+\det \E'= \Ol_{\proj^4}$. It follows that 
$-K_{X'}= 2\xi_{\E'}$ is ample, therefore 
$X'$ is a Fano manifold and $\E'$ is a rank two Fano bundle on $\proj^4$, which,
by \cite[Main Theorem]{APW}, is decomposable as
$\E' \simeq \Ol_{\proj^4}(a) \oplus \Ol_{\proj^4}(b)$.
We can thereby write $\E \simeq (aH -E) \oplus (bH -E)$.
Now, recalling that $\E$ is ample and that $K_Y+\det \E=\Ol_Y$, it is easy to prove that
$(a,b)=(2,3)$.\qed\par

\begin{proposition}\label{XbirYnot} Let $(X,Y,\E)$ be a ruled Fano fivefold of index two 
with $\rho_X = 3$ such that $X$ has a birational contraction but $Y$ has not.
Then one of the following happens:
\begin{enumerate}
\item $X$ is the blow up of a cone in $\proj^9$ over the Segre embedding 
$\pd \times \pd \subset \proj^8$ along its vertex;
\item $X$ is the blow up of $\proj^5$ in two non meeting planes;
\item $X$ is the blow up of a general member of $\Ol(1,1) \subset \pd \times \proj^4$ along
a two dimensional fiber of the second projection.
\end{enumerate}
In these cases the corresponding pairs $(Y,\E)$ are, respectively,
\begin{enumerate}
\item $(\pd \times \pd, \Ol(1,1) \oplus \Ol(2,2))$;
\item $(\pd \times \pd, \Ol(1,2) \oplus \Ol(2,1))$;
\item $(\proj_\pd(T\pd(-1)\oplus \Ol_\pd) \subset \pd \times \pt, \Ol(1,1) \oplus \Ol(1,2))$.
\end{enumerate}
\end{proposition}

{\bf Proof.} \quad First of all it is clear that $X$ cannot be a product
$\pu \times Y$; by proposition \ref{noprod}, recalling that $Y$ has not
birational contractions, the only possible cases are 
$Y \simeq \proj_\pd(T\pd(-1)\oplus \Ol_\pd)$ or $Y \simeq \pd \times \pd$.\par
\smallskip
Let $\vartheta \subset \cone(X)$ be an extremal ray associated to
a birational contraction $\pst:X \to X'$ and let $\theta \subset \cone(Y)$ be its
fellow ray, with associated contraction $\ft:Y \to W$.\\
Denote by $E$ the exceptional locus of $\pst:X \to X'$;
if $E \cdot R_\E =0$, then $E=p^*E_Y$ with $E_Y$ 
an effective divisor on $Y$.
Being $E$ not nef, also $E_Y$ is not nef, and $Y$ has a birational contraction,
against the assumptions. Therefore $E \cdot R_\E >0$ and $E$ dominates $Y$.\\
The fibers of $\pst$ have dimension $\ge 2$ by proposition \ref{fiberlocus}; then,
by lemma \ref{fiberdim}, also the fibers of $\ft$ have dimension $\ge 2$,
hence $\ft$ is a $\pd$-bundle contraction onto $W \simeq \pd$.
By proposition \ref{spcontr2}, $\pst$ is the blow up of a smooth surface $S \subset X'$
and, denoted by $f$ a fiber of $p$, we have $E \cdot f=1$.\\
Let $y$ be a point in $Y$ and let $F_y \simeq \pd$ be the fiber of
$\ft$ through $y$; by the proof of proposition \ref{spcontr2}, 
$\E_{|F_y} \simeq \Ol_\pd(1) \oplus \Ol_\pd(2)$ and $E \cap p^{-1}(F_y)$
is the section corresponding to the $\Ol_\pd(1)$ summand. 
In particular the divisor $E$ cannot contain $f=p^{-1}(y)$. 
It follows that $E$ is a section of $p$, thus $E \simeq Y$.\par
\medskip
Suppose that $X'$ is not a Fano manifold; by \cite[Proposition 3.4]{Wi1}, 
$E$ is negative on another extremal ray
$\ovth \subset \cone(X)$, hence the exceptional locus of the associated 
contraction $\psot:X \to X''$ is contained in $E$, whence $\psot$ is birational.\\
Arguing as above, $\psot:X \to X''$ is the blow up of a smooth fivefold
along a smooth surface, thus its exceptional locus is the divisor $E$;
consequently $E$ has two $\pd$-bundle structures over smooth surfaces and we have 
$E \simeq Y \simeq \pd \times \pd$.\par
\smallskip
Since $E$ is a section of $p$, there exists an exact sequence
$$\shse{\Ol(a_1,a_2)}{\E}{\Ol(b_1,b_2)}$$
such that $E \simeq \xi_\E \otimes p^*\Ol(-a_1,-a_2)$; being 
$E \cdot \vartheta = E \cdot \ovth = -1$, we have
$a_1=a_2=2$; then 
$$-1=E \cdot \vartheta = (1-a_1) =E \cdot \ovth = (1-a_2).$$
Recalling that $\det \E = -K_{\pd \times \pd}= \Ol(3,3)$, we obtain $b_1=b_2=1$;
since $h^1(\pd \times \pd, \Ol(a_1-b_1,a_2-b_2))= h^1(\pd \times \pd, \Ol(1,1))=0$,
the above sequences splits, the vector bundle $\E$ is decomposable: 
$\E \simeq \Ol(1,1) \oplus \Ol(2,2)$, and we are in case (1).\par
\medskip
We can now assume that $X'$ is a Fano manifold; consider the commutative diagram
as in \ref{gc}
$$
\xymatrix@=40pt{X  \ar[r]^{\pst} \ar[d]_{p}  &X' \ar[d]^{p'}\\
Y \ar[r]_\ft  &\pd}
$$

Let $x \in \pd$ be a general point; the fibers $G=p'^{-1}(x)$ and $F=\pst^{-1}(p'^{-1}(x))$
are smooth and, by the commutativity of the diagram, 
$F=p^{-1}(\ft^{-1}(x))= \proj_\pd(\Ol_\pd(1)\oplus \Ol_\pd(2))$; therefore $G \simeq \pt$.\\
By lemma \ref{brauer} there exists a rank four vector bundle $\F$ over 
$\pd$ such that $X'=\proj_\pd(\F)$; in particular $\F$ is a Fano bundle over $\pd$.\\
By the canonical bundle formula for blow ups we have
$$-\pst^*K_{X'}=-K_X +2E=2(\xi_\E+E),$$
whence the index of $X'$ is two. Writing $K_{X'}$ with the canonical bundle formula
for projectivizations
$$K_{X'}=-4\xi_\F +p'^*(\Ol_\pd(-3)+c_1(\F)),$$
this implies that the first Chern class of $\F$ is odd. 
By the classification in \cite{SW1} either 
$\F \simeq \Ol_\pd^{\oplus 3} \oplus \Ol_\pd(1)$ or 
$\F \simeq T\pd(-1)\oplus \Ol^{\oplus 2}_\pd$.\par
\smallskip
As for every $x \in \pd$ the fiber $F_x=\pst^{-1}(p'^{-1}(x))$ is the blow up of
$\pt$ at a point and the fiber $G_x=p'^{-1}(x)$ is a projective space of dimension three,
we have that $S$, the center of the blow-up $\pst$, is a section of $p'$; therefore
we have an exact sequence
\begin{equation}\label{centerseq}
\shse{\G}{\F}{\Ol(a)}
\end{equation}
such that $S$ is the zero locus of a section of the vector bundle 
$\xi_\F \otimes p'^*{\G}^\vee$; in particular the conormal bundle $N^*_{S/X'}$ of 
$S$ is $(p'^*\G \otimes \xi_\F^{-1})_{|S}$. Recall that the exceptional divisor $E$ is 
the projectivization of the conormal bundle of
$S$, i.e. $E \simeq \proj_S(N^*_{S/X})$.\par
\medskip
If $E \simeq Y \simeq \pd \times \pd$, then 
$N^*_{S/X}$,  hence $\G$ is decomposable.
It follows that $h^1(\G(-a))=0$, thus the sequence splits
and we have $\G \simeq \Ol_\pd^{\oplus 3}$, $\F \simeq \Ol_\pd^{\oplus 3} \oplus \Ol_\pd(1)$,
i.e. $S$ is the section corresponding to
the surjection $\F \to \Ol_\pd(1)$ and it is disjoint from the exceptional
divisor of the blow down $X' \to \proj^5$. We thereby conclude that $X$
is the blow up of $\proj^5$ in two non meeting planes.\par
\medskip
Suppose now that $E \simeq Y \simeq \proj_\pd(T\pd(-1)\oplus \Ol_\pd)$.\\
Let $\oth$ be the extremal ray corresponding to the contraction
$\fot:Y \to \pt$, which is a special B\v anic\v a scroll, and let
$\psot:X \to Z$ be the contraction associated to $\ovth$, the fellow ray
of $\oth$; by proposition \ref{spcontr3} $\psot$ is a special B\v anic\v a scroll.\\
Let $\sigma \subset \cone(X)$ be the face spanned by $\vartheta$ and $\ovth$;
the contraction of this face, call it $\psi_\sigma$, factors through the contraction
$\pst:X \to X'$ and we have a commutative diagram
$$
\xymatrix@=20pt{
\pt  & & Z \ar[rd] &\\
 Y \ar[d]_{\ft} \ar[u]^{\fot} &  X \ar[l]_{p} \ar[rr]^{\psi_\sigma}
 \ar[ru]^{\psot} \ar[rd]_{\pst} & & W'\\
  \pd  & & X' \ar[ll]^{p'} \ar[ru]_{\pi}  & }
$$
The morphism $\pi:X' \to W'$ is the contraction of $X'$ different from the projection
onto $\pd$; since $\dim W' \le \dim Z < \dim X$, $\pi$ is a fiber type contraction,
so $\F \simeq T\pd(-1) \oplus \Ol_\pd^{\oplus 2}$ and $W' \simeq \pt$.\\
We claim that $E \cdot \ovth =0$; indeed, if this is not true,
then, for every $x \in X$, denoting by $(F_{\ovth})_x$ the fiber of
$\psot$ containing $x$, we will have
$$\dim \psot^{-1}(\psot(\pst^{-1}(\pst((F_\ovth)_x)))) \ge 3.$$

Denoting by $V^\vartheta$ and $V^\ovth$ the families
of minimal degree rational curves whose numerical class is
in $\vartheta$ and $\ovth$, respectively, and by $(F_\sigma)_x$ the fiber
of $\psi_\sigma$ containing $x$ we will have
$$(F_\sigma)_x \supset \cloc(V^\vartheta,V^\ovth)_x \supset 
\psot^{-1}(\psot(\pst^{-1}(\pst((F_\ovth)_x)))),$$
a contradiction, since the general fiber of $\psi_\sigma:X \to \pt$
is two dimensional.\par
\medskip
As we have already noticed, $E=\proj(N^*_{S/X'})$ and, since $E \simeq Y$, 
$N^*_{S/X'} \simeq T\pd(b-1) \oplus \Ol_\pd(b)$ for some $b$.
The fact that $E \cdot \ovth=0$ implies that $b=0$, so 
$$\G \simeq (p'^*\G)_{|S} \simeq (\xi_\F)_{|S} \otimes N^*_{S/X'} \simeq (\xi_\F)_{|S} 
\otimes(T\pd(-1)\oplus \Ol_\pd) \simeq T\pd(x-1)\oplus \Ol_\pd(x)$$
with $x \ge 0$ since $\F$ is nef; by the sequence (\ref{centerseq}),
we have an injection
$$0 \to T\pd(x-1)\oplus \Ol_\pd(x) \to \Ol_\pd(-1)\oplus \Ol_\pd^{\oplus 2},$$
which forces $x =0$. It follows that $S$ corresponds to a surjection
$\F \to \Ol_\pd \to 0$, so it is a two dimensional fiber of the special  B\v anic\v a 
scroll contraction of $X'$.\qed\par


\section{Proof of theorem \ref{rho2}}

The main idea of the proof of theorem \ref{rho2} is to consider,
when possible, a smooth divisor $Y'$ in the linear system of the ample generator
of $Y$, and to study the manifold $X'=\proj_{Y'}(\E_{|Y'})$;
in order to do that we first establish some relations between
the geometry of $X$ and the geometry of $X'$.

\begin{lemma}\label{dpnef}
Let $Y$ be a smooth variety, $L \in \pic(Y)$ an ample line bundle
and $Y' \in |L|$ an effective divisor.
Let $\E$ be a rank two vector bundle on $Y$ and denote by $\E_{Y'}$
its restriction to $Y'$. Then
\begin{enumerate}
\item[a)] if $\E_{Y'}$ is spanned, then $\E$ is nef; 
\item[b)] if $h^i(\E_{Y'}(-jL))=0$ for  $i=0,1$ and every $j \ge 1$,
then $H^0(Y,\E) \simeq H^0(Y',\E_{Y'})$.
\end{enumerate}
\end{lemma}

{\bf Proof.} \quad 
By definition, the nefness of $\E$ is the nefness of its tautological bundle;
let $X=\proj_Y(\E)$ and let $X'=\proj_{Y'}(\E_{Y'})$.
Since the restriction of $\xi_\E$ to $X'$ is spanned, if $\xi_\E \cdot C <0$
for some effective curve $C$, then $C \cap X' = \emptyset$.
By the ampleness of $Y'$ in $Y$ this implies that $C$ is a fiber
of the natural projection $p:X \to Y$, but this is impossible
since such curves cover $X$.\\
To prove b), by the exact sequence
$$\shse{\E(-L)}{\E}{\E_{{Y'}}},$$
we have to show that $h^0(\E(-L))= h^1(\E(-L))=0$, and
this follows from \cite[Corollary 4.1.6]{BS2}.\qed\par

\begin{proposition}\label{dpcones}
Let $Y$ be a smooth variety of Picard number one and dimension $\ge 4$, 
$\E$ a rank two vector bundle on $Y$,
$L \in \pic(Y)$ an ample line bundle and $Y' \in |L|$ an effective divisor.\\
Assume that $\E_{Y'}=\E_{|Y'}$ is spanned and that $|\xi_{\E_{Y'}}|$
defines an extremal contraction $\f_{\vartheta'}:X'=\proj_{Y'}(\E_{Y'}) \to Z$ associated
to an extremal ray $\vartheta' \subset \cone(X')$.
Then, under the identification $N_1(X') \simeq N_1(X)$, given by the inclusion
$i:X' \to X$, we have $\cone(X')=\cone(X)$.
\end{proposition}

{\bf Proof.} \quad Since $\dim Y \ge 4$, by Weak Lefschetz theorem we have $\rho_{Y'}=1$,
hence the cones of curves $\cone(X)$ and $\cone(X')$ have dimension
two and, under the identification $N_1(X') \simeq N_1(X)$, they have in common
the extremal ray $R_\E$ corresponding to the bundle projection.
We have therefore to prove $\vartheta'$ is extremal in $\cone(X)$, too.\\
Since $(\xi_\E)_{|X'}=\xi_{\E_{Y'}}$ is zero on $\vartheta'$, 
if $\vartheta'$ is not extremal in $\cone(X)$
we have $\xi_\E \cdot C <0$ for some curve whose class is in
$\cone(X) \setminus \cone(X')$. This contradicts the fact that,
by lemma \ref{dpnef} a), $\E$ has to be nef.\qed\par

\begin{corollary}\label{extcontr} Let $(X,Y,\E)$ be a ruled Fano
fivefold of index two and Picard number $\rho_X=2$, let $L$ be the
ample generator of $\pic(Y)$, and assume that there exists an effective divisor
$Y' \in |L|$ such that $\E_{Y'}=\E_{|Y'}$ is spanned and that $|\xi_{\E_{Y'}}|$
defines an extremal contraction $\f_{\vartheta'}:X \to Z$ of fiber type.
Then there exists an extremal contraction 
$\psi_\vartheta:X \to Z$ such that $(\psi_{\vartheta})_{|X'}=\f_{\vartheta'}$.
\end{corollary}

{\bf Proof.} \quad This assertion follows from \cite[Proposition 3.13]{AO1}.\qed \par
\bigskip

{\bf Proof of theorem \ref{rho2}.} \quad By lemma \ref{bundlebase}, $Y$ is a Fano variety of
pseudoindex $i_Y \ge i_X = 2$; moreover, since $\rho_X=2$, we have $\rho_Y=1$.\par
\medskip
If $r_Y=i_Y=2$, i.e. $Y$ is a Mukai manifold, then, denoted by $\Ol_Y(1)$
the ample generator of $\pic(Y)$, by \cite[Theorem 1]{Me} 
a general section $Y'$ in $|\Ol_Y(1)|$ is smooth, and so
it is a Fano threefold of index one. 
By adjunction $X'=\proj_{Y'}(\E_{Y'})$ is a Fano manifold,
hence we can apply \cite[Theorem 8.4]{Lan} to get $X'= \pu \times  Y'$.\\
Up to a twist, we can assume that $\E_{Y'} \simeq \Ol_{Y'} \oplus \Ol_{Y'}$; 
this bundle verifies the assumptions of proposition \ref{dpcones}, so, by
corollary \ref{extcontr}, there exists
an extremal contraction $\psi_\vartheta:X \to \pu$; by lemma \ref{products} we have
$X \simeq \pu \times Y$.\par
\smallskip
If $r_Y=i_Y=3$, i.e. $Y$ is a del Pezzo manifold, we again denote by $\Ol_Y(1)$
the ample generator of $\pic(Y)$ and we take a general divisor 
$Y' \in |\Ol_Y(1)|$. By adjunction $X'=\proj_{Y'}(\E_{Y'})$ is a Fano manifold;
by \cite[Theorem 8.2]{Lan} and \cite[Proposition 4.2]{SW4} we have the 
following possibilities for $(Y', \E_{Y'})$ (here the vector bundles are not
normalized as in definition \ref{rfmi}):
\begin{enumerate}
\item $(V_d,\Ol_{V_d} \oplus \Ol_{V_d}(-1))$;
\item $(V_4,$ restriction of a spinor bundle on $\mathbb Q^4)$;
\item $(V_5,$ restriction of the universal bundle on $G(1,4))$.
\end{enumerate}\par

\smallskip
{\bf Case 1} \quad $(Y', \E_{Y'}) \simeq (V_d,\Ol_{V_d} \oplus \Ol_{V_d}(-1))$.\par
\medskip
By lemma \ref{dpnef} b) $H^0(Y,\E)\simeq H^0(Y',\E) \simeq \mathbb C$.
It follows that $\E$ has a section, $s$; this section does not vanish on ${Y'}$, which
is ample, whence $s$ can vanish only at points outside ${Y'}$. Let $x$ be one of these points
and let $l$ be a line through $x$; $\E(1)$ is ample and $\det \E(1) \simeq \Ol_Y(3)$, 
so $\E$ restricts to $l$ as $\Ol_{\pu} \oplus \Ol_\pu(-1)$, and
$s$ cannot vanish on $l$.\\
We thereby have a short exact sequence
$$\shse{\Ol}{\E}{L}$$
where, computing the splitting type, we have $L=\Ol_Y(-1)$; consequently the sequence splits
and $\E \simeq \Ol_Y \oplus \Ol_Y(1)$.\par
\medskip
{\bf Case 2} \quad $(Y', \E_{Y'}) \simeq (V_4,$ restriction of 
a spinor bundle on $\mathbb Q^4)$.\par
\medskip
In case (2), as proved in \cite[4.4]{SW4}, $X'$ has a conic bundle structure
$\f:X' \to \pt$, and can be described as a divisor in the flag manifold of lines and
points in $G(1,3) \times \pt$.
Indeed, $\f_*\xi_{\E_{Y'}(1)} \simeq \Omega \pt(3)$ and
the flag manifold can be identified with the projectivization $\proj_\pt(\Omega \pt(3))$;
with this description $X'$ is a divisor in $|2\xi_{\Omega\pt(3)}-2\f^*\Ol_\pt(1)|$.\par
\smallskip
Since $\E$ is spanned on ${Y'}$ and $|\xi_{\E_{Y'}}|$ defines a fiber
type contraction, by corollary \ref{extcontr},
there exists a contraction $\pst:X \to \pt$ such that its restriction
to $X'$ is the conic bundle contraction $\f:X' \to \pt$.\\
In particular, since the restriction of $\pst$ to $X'$ is equidimensional
and $X'$ is $\pst$-ample, also $\pst$ is equidimensional and,
by adjunction, is a quadric bundle contraction.\par
\smallskip
Let $\F={\pst}_*\xi_{\E(1)}$; $\F$ is a vector bundle of rank four and $X$ embeds in 
$\proj_{\proj^3}(\F)$ as a divisor of relative degree $2$, i.e.
$X \in |2\xi_{\E(1)} +\pst^*\Ol_\pt(x)|$.\\
The vector bundle $\F$ has $\G =\f_*\xi_{\E_{Y'}(1)} \simeq \Omega \pt(3)$
as a quotient. Indeed, if $x \in \proj^3$ is a point and we denote by $F$ and $f$
the fibers of $\pst$ and ${\pst}_{|X'}=\f$ over $x$, we have
that $\G_x =H^0((\xi_\E)_{|f})$ is a quotient of $\F_x= H^0((\xi_\E)_{|F})$.\\
It follows that there exists an exact sequence on $\pt$:
$$\shse{\Ol(a)}{\F}{\Omega\pt(3)}.$$

Since $(\xi_{\E(1)})_{|X'}=\xi_{\Omega\pt(3)}$, $(\pst^*\Ol(1))_{|X'}=\f_\pt^*\Ol_{\pt}(1)$
and $X_{|\proj(\G)}= X'$, we have $x=-2$ and
$X \in |2\xi_{\E(1)} -2\psi^*\Ol_\pt(1)|=|2\xi_\E|$.
By adjunction
$$-2\xi_{\E(1)} =K_X = (K_{\proj_\pt(\F)}+X)_X= -2\xi_{\E(1)} +\psi^*\Ol_\pt(c_1(\F)-6),$$
hence $c_1(\F)=6$. Computing the degree in the above sequence, we have $a=1$. Therefore
the sequence splits and we have 
$\F \simeq \Omega\proj^3(3) \oplus \Ol_\pt(1)$.\par
\medskip
{\bf Case 3}\quad $(Y', \E_{Y'}) \simeq (V_5,$ restriction of the universal 
bundle on $G(1,4))$.\par
\medskip

We claim that $\E$ is spanned on $Y$; to prove the claim we show that $\xi_{\E}$ is spanned on 
$X=\proj(\E)$.
Assume that $\bar x \in X$ is a base point of $|\xi_\E|$; since
$\Ol_Y(1)$ is very ample, we can find a smooth section $Y'' \in |\Ol_Y(1)|$ 
containing $p(x)$.
The restriction 
${(\xi_{\E})}_{|Y''}=\xi_{\E_{Y''}}$ is spanned, so there exists a section 
of $(\xi_{\E})_{Y''}$ which does not vanish at $x$ and this section, 
by lemma \ref{dpnef} b), extends to $X$.\\
We have thus proved that $\E$ is spanned; again by lemma \ref{dpnef} b), 
$h^0(Y,\E)=h^0(Y'', \E_{Y''})=5$ so 
we have an exact sequence of vector bundles
$$\shse{{\mathcal G}}{\Ol_Y^{\oplus 5}}{\E}$$
which gives an injection $X \to \proj^4 \times Y$ and then an injection 
$X \to \proj^4 \times G(1,4)$. We claim that $X$ is the intersection of 
$p^{-1}(Y)$ with the flag manifold of lines and points
in $G(1,4) \times \proj^4$.
Indeed, given a point $y \in Y$, denoting by $Y'$ a smooth member of $\Ol_Y(1)$
passing through $y$, $\E_{Y'}$ is the restriction of the universal bundle of
$G(1,4)$, thus the fiber of $\E$ over $y$ is the line parametrized by $y \in G(1,4)$.\par
\smallskip
If $i_Y=4$ then, by \cite[Theorem 0.1]{Mi2}, $Y \simeq \quadr^4$. We can apply
\cite[Theorem 2.4]{APW} to get that  $\E$ is decomposable (the other bundles
have odd $c_1$, while, in our case, since $\det \E= -K_Y = \Ol_{\mathbb Q^4}(4)$,
$c_1(\E)$ is even) and we are in case (2) of theorem \ref{rho2}.\par
\smallskip
If $i_Y=5$ then, by \cite[Corollary 0.4]{CMS} or \cite[Theorem 1.1]{Ke}, 
$Y \simeq \proj^4$. We can apply \cite[Theorem 2.4]{APW} to get that $\E$ is decomposable, 
hence we are in case (1) of theorem \ref{rho2}. Note that only the bundles whose
projectivization gives a Fano manifold of index two are considered.
\qed \par

\providecommand{\bysame}{\leavevmode\hbox to3em{\hrulefill}\thinspace}

\end{document}